\documentclass[review,12pt]{elsarticle}

\usepackage{amssymb}
\usepackage{amsmath}
\usepackage{subfigure}
\usepackage{color}
\usepackage{mathrsfs} 
\usepackage[colorlinks=true]{hyperref}
\usepackage{amsthm}

\usepackage{bm}
\usepackage{cleveref}
\usepackage{multirow}

\usepackage{algorithm}
\usepackage{algorithmicx}
\usepackage{algpseudocode}

\usepackage{booktabs}
\usepackage[normalem]{ulem}

 \newtheorem{thm}{Theorem}[section]

 \theoremstyle{definition}
 
 \theoremstyle{remark}
 \newtheorem{rem}[thm]{Remark}

\def\ep{\epsilon}
\def\pt{\partial}

\def\hf{\frac{1}{2}}

\def\Re{\mathrm{Re}}
\def\Ma{\mathrm{Ma}}
\def\Pr{\mathrm{Pr}}

\def\bfu{\mathbf{u}}
\def\bfv{\mathbf{v}}
\def\bfw{\mathbf{w}}

\def\bfx{\mathbf{x}}
\def\bff{\mathbf{f}}
\def\bfg{\mathbf{g}}

\def\bfB{\mathbf{B}}

\def\bfH{\mathbf{H}}
\def\bfQ{\mathbf{Q}}

\def\bfbu{\mathbf{\bar{u}}}
\def\bfbsig{\bm{\bar{\sigma}}}
\def\bfbvpi{\bm{\bar{\varpi}}}
\def\bftQ{\mathbf{\widetilde{Q}}}

\def\vb{\bar{v}}

\def\bftf{\widetilde{\bff}}
\def\bftg{\widetilde{\bfg}}

\def\vu{\vec{u}}
\def\vv{\vec{v}}
\def\vT{\vec{T}}

\def\vb{\vec{b}}
\def\fM{\mathcal{M}}

\begin{document}

\begin{frontmatter}



\title{A Second-Order Relaxation Flux Solver for Compressible Navier-Stokes Equations based on Generalized Riemann Problem Method }

\author[HKUST]{Tuowei Chen}

\ead{tuowei_chen@163.com}
\cortext[cor1]{Corresponding Author}

\author[IAPCM]{Zhifang Du\corref{cor1}}
\ead{du_zhifang@iapcm.ac.cn}

\address[HKUST]{Department of Mathematics, Hong Kong University of Science and Technology, \\ Clear Water Bay, Kowloon, HongKong}

\address[IAPCM]{Institute of Applied Physics and Computational Mathematics, \\ 100048 Beijing, PR China}

\begin{abstract}
In the finite volume framework, a Lax-Wendrof type second-order flux solver for the compressible Navier-Stokes equations is proposed by utilizing a hyperbolic relaxation model. 
The flux solver is developed by applying the generalized Riemann problem (GRP) method to the relaxation model that approximates the compressible Navier-Stokes equations.
The GRP-based flux solver includes the effects of source terms in numerical fluxes and treats the stiff source terms implicitly, allowing a CFL condition conventionally used for the Euler equations. The trade-off is to solve linear systems of algebraic equations. The resulting numerical scheme achieves second-order accuracy within a single stage, and the linear systems are solved only once in a time step. The parameters to establish the relaxation model are allowed to be locally determined at each cell interface, improving the adaptability to diverse flow regions. Numerical tests with a wide range of flow problems, from nearly incompressible to supersonic flows with strong shocks, for both inviscid and viscous problems, demonstrate the high resolution of the current second-order scheme.
\end{abstract}




\end{frontmatter}


\vspace{2mm}
\section{Introduction}\label{sec:intro}
In the application of computational fluid dynamics (CFD), the compressible Navier-Stokes equations serve as the governing partial differential equations (PDEs) of compressible viscous flows. Over the past decades, finite volume schemes commonly used for the compressible Navier-Stokes equations primarily employ the method of lines, where the spatial discretizations for the convective and viscous fluxes are treated separately \cite{Blazek2015}. Typically, the convective flux is discretized based on an exact or approximate Riemann solver, while the viscous flux is treated with a central finite difference method. The time integration is then implemented using the multi-stage Runge-Kutta (RK) methods, which are widely used for their simplicity. However, a fully explicit numerical scheme suffers from a severe time step restriction due to diffusion terms, particularly when the mesh size varies dramatically. In contrast, a fully implicit scheme could be computationally expensive since it leads to large nonlinear systems to solve. A typical strategy to address this issue is the implicit-explicit (IMEX) method \cite{BBMRV2015,GMPT2021,RGMC2013,SSHM2010}. Developed based on Riemann solvers, RK type IMEX schemes involve solving systems of algebraic equations at each intermediate stage \cite{ARS1997,CDN2001,WWZS2016}.

An alternative approach is provided by Lax-Wendroff (LW) type flux solvers \cite{LW1960}. 
LW type solvers are numerical implementations of the Cauchy-Kowalevski approach in the context of PDEs. Apart from Riemann solutions, it also obtains time derivatives of flow variables at cell interfaces, allowing a high-order interpolation of fluxes in the time direction. This is done by converting space derivatives of flow variables into time ones. Corresponding numerical schemes achieve high-order time accuracy through a single-stage evolution. Among the LW type flux solvers for hyperbolic PDEs, the generalized Riemann problem (GRP) solver \cite{BF1984,BF2003,BLW2006} is a typical representative. The GRP is an initial value problem (IVP) of the governing PDEs considering piecewise smooth initial data. 
While LW type flux solvers are extensively applied to hyperbolic PDEs, their application to solving convection-diffusion PDEs is very limited. A successful numerical scheme for the compressible Navier-Stokes equations developed within the LW framework is the gas-kinetic scheme (GKS) \cite{Xu2001}, which is based on the kinetic Bhatnagar-Gross-Krook (BGK) model. 
Developing the GRP solvers for the compressible Navier-Stokes equations without utilizing kinetic models remains a significant challenge, as there is currently no established method to directly resolve the time evolution of viscous fluxes.
To address this, in this paper we turn to solving the GRP of a hyperbolic relaxation model that approximates the original compressible Navier-Stokes equations.

In our previous work \cite{CL2023}, we extended the Jin-Xin hyperbolic relaxation model for inviscid conservation laws \cite{JX1995} to viscous ones. Using the resulting models, we also investigated the GRP-based second-order finite volume schemes for solving convection-diffusion PDEs.
However, this approach encountered two significant challenges. First, as with other relaxation-based schemes \cite{Nishikawa2011,Nishikawa2014,LN2017,LLLN2019,MT2014,DPRZ2016,WA2024a,WA2024b}, updating both conservative and relaxation variables requires solving additional PDEs and storing extra data. Second, the numerical scheme faces a severe time step restriction of $\Delta t \sim O(\sqrt{\epsilon}\Delta x)$, where $\epsilon$ is the relaxation parameter. These limitations motivate the present work, which seeks to overcome flux stiffness and avoid finite volume discretizations of relaxation variable equations.

This paper aims to develop a LW type flux solver and the corresponding finite volume scheme for the compressible Navier-Stokes equations.
To achieve this, a hyperbolic relaxation model with stiff components regarded as source terms is employed. 
The flux solver for the compressible Navier-Stokes equations is developed by solving the GRP of the hyperbolic relaxation model. In the procedure of computing time derivatives, this flux solver handles the stiff source terms in the relaxation model implicitly.
From the viewpoint of the relaxation model, this implicit treatment effectively addresses the stiffness of the model. As a trade-off, linear systems of algebraic equations are solved at each time step. 
From the viewpoint of the original compressible Navier-Stokes equations, solving linear systems acts as the implicit discretization to the diffusion terms in RK type IMEX schemes, thereby overcoming the parabolic time step restriction.
The coefficient matrix is sparse and strictly diagonally dominant, making the system efficient to solve.

The resulting numerical scheme naturally inherits the advantages of GRP-based schemes. Compared with the Riemann-solver-based RK type IMEX schemes, it achieves second-order accuracy in a single stage. As a consequence, a higher efficiency is achieved since the data reconstruction step is reduced and the linear systems are solved only once at each time step. Moreover, the present scheme is more compact compared with multi-stage time-stepping schemes, improving the resolution to multi-dimensional structures in the flow field.

Compared with previous relaxation schemes \cite{Nishikawa2011,Nishikawa2014,LN2017,LLLN2019,MT2014,DPRZ2016,WA2024a,WA2024b}, 
a key distinguishing feature of the proposed method is its innovative utilization of relaxation variables: they are solely employed for constructing second-order accurate numerical fluxes, and the necessity to update these relaxation variables in numerical discretizations is avoided.
This feature makes it operate as a flux solver rather than a full numerical scheme. This consequently reduces both the computational cost and the memory consumption.
Meanwhile, the parameters involved in the relaxation model can vary at different cell interfaces and at each time step. Rather than being a fixed small constant, the relaxation parameter $\ep$ takes a value considering the discontinuity of flow variables on the two sides of the cell interface. In the smooth region, $\ep \ll \Delta t$, and the relaxation effect is dominant in the model.
In the presence of discontinuities, $\ep$ increases, and the convective effect supplies numerical dissipation to suppress non-physical oscillations. Numerical experiments show that this numerical treatment enhances the robustness of the scheme. This numerical procedure that produces numerical fluxes of second-order time accuracy for the compressible Navier-Stokes equations is termed as the \emph{Relaxation Flux Solver (RFS)} in the present paper.

The rest of this paper is organized as follows. In \Cref{preliminary}, a hyperbolic relaxation model for the two-dimensional (2-D) compressible Naiver-Stokes equations is briefly reviewed. In \Cref{FVCNS}, a GRP-based finite volume scheme is proposed, where the relaxation model is employed to compute the numerical flux. In \Cref{GRPFS}, the GRP solver for the relaxation model is developed, which can be regarded as a LW type flux solver for the 2-D compressible Naiver-Stokes equations. The numerical scheme together with the flux solver is summarized as an algorithm in \Cref{sec:algorithm}. Numerical results are presented in \Cref{Examples} to demonstrate the high-resolution performance of the current scheme. The last section presents discussions.

\section{Governing equations and relaxation model} \label{preliminary}
The 2-D compressible Navier-Stokes equations can be expressed as
\begin{equation}
	\frac{\pt\bfu}{\pt t}+ \frac{\pt \bff_c(\bfu)}{\pt x}+\frac{\pt\bfg_c(\bfu)}{\pt y}
	=\frac{\pt \bff_v(\bfu,\nabla \bfu)}{\pt x}+\frac{\pt \bfg_v(\bfu,\nabla \bfu)}{\pt y},
	\label{CNS0}
\end{equation}
where $\bfu=(\rho,\rho u,\rho v, \rho E)^\top$ is the vector of conservative variables, $\rho$ is the density, $u$ and $v$ are the velocity components in $x$- and $y$-directions, and $E$ is the specific total energy. The convective and viscous fluxes are defined as 
\begin{equation}
	\label{flux}
	\begin{array}{l}
		\bff_c(\bfu)=\begin{pmatrix}
			\rho u, \ \rho u^2+p, \rho uv, \ (\rho E+p)u
		\end{pmatrix}^\top,	\\
		\bfg_c(\bfu)=\begin{pmatrix}
			\rho v, \ \rho uv, \ \rho v^2+p, \ (\rho E+p)v
		\end{pmatrix}^\top, \\
		\bff_v(\bfu,\nabla\bfu)=\begin{pmatrix}
			0, \ \tau_{11}, \ \tau_{12}, \ \tau_{11}u+\tau_{12}v+q_1
		\end{pmatrix}^\top,\\
		\bfg_v(\bfu,\nabla\bfu)=\begin{pmatrix}
			0, \ \tau_{21}, \ \tau_{22}, \ \tau_{21}u+\tau_{22}v+q_2
		\end{pmatrix}^\top.
	\end{array}
\end{equation}
The pressure $p$ is determined by the ideal gas equation of state (EOS):
\begin{equation}\label{eq:EOS}
	p=(\gamma-1)\big[\rho E-\frac{1}{2}\rho(u^2+v^2)\big],
\end{equation} 
where $\gamma$ is the specific heat ratio.
For Newtonian fluids, based on the Stokes hypothesis and the Fourier law, the viscous stress and the heat flux are given by
\begin{equation*}
	\begin{aligned}
		&\tau_{11}=\frac{4}{3}\mu\frac{\pt u}{\pt x}
		-\frac{2}{3}\mu\frac{\pt v}{\pt y},
		\quad \tau_{12}=\tau_{21}
		=\mu\bigg(\frac{\pt u}{\pt y}+\frac{\pt v}{\pt x}\bigg),\quad
		\tau_{22}=-\frac{2}{3}\mu\frac{\pt u}{\pt x}
		+\frac{4}{3}\mu\frac{\pt v}{\pt y},\\
		&q_1=\kappa\frac{\pt T}{\pt x},\quad q_2=\kappa\frac{\pt T}{\pt y}, 
		\quad \kappa = \frac{\gamma\mu}{\Pr(\gamma-1)}.
	\end{aligned}
	\label{visflux1}
\end{equation*}
In the above equations, $T$ is the temperature, $\mu$ is the dynamic viscosity, $\kappa$ is the heat conductivity, and $Pr$ is the Prandtl number. By the ideal gas EOS \eqref{eq:EOS}, the temperature $T$ is defined as
\begin{equation}\notag
	T=\frac{p}{\rho}.
\end{equation}
The viscosity coefficient $\mu$ generally depends on the temperature $T$ and can take any reasonable form. 

The relaxation model considered in \cite{CL2023} is written as 
\begin{equation}	\label{CNS1}
	\left.\begin{cases}
		\begin{aligned}
			&\frac{\pt\bfu}{\pt t}+\frac{\pt \bfv}{\pt x}+\frac{\pt\bfw}{\pt y}=0, \\
			&\frac{\pt\bfv}{\pt t}+\bigg[\frac{1}{\ep}\bfB_{11}(\bfu)+a^2\bigg] \frac{\pt \bfu}{\pt x}
			+\frac{1}{\ep}\bfB_{12}(\bfu) \frac{\pt \bfu}{\pt y}
			=\frac{1}{\epsilon}\big(\bff_c(\bfu)-\bfv\big),\\
			&\frac{\pt\bfw}{\pt t}+\frac{1}{\ep}\bfB_{21}(\bfu) \frac{\pt \bfu}{\pt x}
			+\bigg[\frac{1}{\ep}\bfB_{22}(\bfu)+a^2\bigg] \frac{\pt \bfu}{\pt y}
			=\frac{1}{\epsilon}\big(\bfg_c(\bfu)-\bfw\big),
		\end{aligned}
	\end{cases}\right.
\end{equation} 
where $a>0$, $\ep>0$ is a small parameter standing for the relaxation time, $\bfv$ and $\bfw$ are relaxation variables for fluxes in the $x$- and $y$-directions, repectively. The matrices $\bfB_{ij}(\bfu)$ $(i,j=1,2)$ (see \ref{tensor} for {detailed expressions) satisfy
\begin{equation} \label{matricesB}
	\begin{aligned}
		&\bff_v(\bfu,\nabla\bfu)
		=\bfB_{11}(\bfu)\frac{\pt\bfu}{\pt x}+\bfB_{12}(\bfu)\frac{\pt\bfu}{\pt y}, \\
		&\bfg_v(\bfu,\nabla\bfu)
		=\bfB_{21}(\bfu)\frac{\pt\bfu}{\pt x}+\bfB_{22}(\bfu)\frac{\pt\bfu}{\pt y}.
	\end{aligned}
\end{equation}
The subcharacteristic condition for the coefficient $a$ has been analyzed in the previous work \cite{CL2023}.

In the present work, we re-formulate the hyperbolic relaxation model \eqref{CNS1} as
\begin{equation}
	\left.\begin{cases}
		\begin{aligned}
			&\frac{\pt\bfu}{\pt t}+\frac{\pt \bfv}{\pt x}+\frac{\pt\bfw}{\pt y}=0, \\
			&\frac{\pt\bfv}{\pt t}+a^2\frac{\pt \bfu}{\pt x}
			=\frac{1}{\epsilon}\bigg[\bff_c(\bfu)-\bff_v(\bfu,\nabla\bfu)-\bfv\bigg],\\
			&\frac{\pt\bfw}{\pt t}+a^2\frac{\pt \bfu}{\pt y}
			=\frac{1}{\epsilon}\bigg[\bfg_c(\bfu)-\bfg_v(\bfu,\nabla\bfu)-\bfw\bigg],
		\end{aligned}
	\end{cases}\right.
	\label{CNS2}
\end{equation}
where $\frac{1}{\ep}\bff_v(\bfu,\nabla\bfu)$ and $\frac{1}{\ep}\bfg_v(\bfu,\nabla\bfu)$ are regarded as stiff source terms.

\section{Finite volume scheme based on hyperbolic relaxation model} \label{FVCNS}
In this section, based on the hyperbolic relaxation model \eqref{CNS2}, a finite volume scheme is proposed for the 2-D compressible Navier-Stokes equations \eqref{CNS0}. The proposed scheme achieves second-order accuracy by discretizing the relaxation model with the GRP methodology.
The definition of the GRP at cell interfaces is provided afterward.

\subsection{Second-order finite volume scheme}\label{sec:fv}
For simplicity, we consider only the case of rectangular meshes in the present paper. The computational cell is $\Omega_{i,j}=[x_{i-\hf}, x_{i+\hf}]\times[y_{j-\hf},y_{j+\hf}]$ with the cell size $\Delta x_i=x_{i+\hf}-x_{i-\hf}$, $i=1,\dots,I$, and $\Delta y_j=y_{j+\hf}-y_{j-\hf}$, $j=1,\dots,J$. The center of $\Omega_{i,j}$ is $(x_i,y_j)$, where $x_i=\frac{1}{2}(x_{i-\hf}+x_{i+\hf})$ and $y_j=\frac{1}{2}(y_{j-\hf}+y_{y+\hf})$. In the interval $[t^n,t^{n+1}]$, the time step is defined as $\Delta t=t^{n+1}-t^n$.

Based on the relaxation model \eqref{CNS2}, a second-order finite volume scheme for the governing equation of $\bfu$ is constructed as
\begin{equation}
	\bfbu^{n+1}_{i,j}=\bfbu^{n}_{i,j}
	-\frac{\Delta t}{\Delta x_i}\big(\bfv^{n+\hf}_{i+\hf,j}-\bfv^{n+\hf}_{i-\hf,j}\big)	-\frac{\Delta t}{\Delta y_j}\big(\bfw^{n+\hf}_{i,j+\hf}-\bfw^{n+\hf}_{i,j-\hf}\big),
	\label{FV}
\end{equation}
where $\bfbu_{i,j}^n$ is the cell average of $\bfu$ over $\Omega_{i,j}$ at $t=t^n$.
Here, the numerical fluxes are approximated by} the mid-point values of $\bfv$ and $\bfw$ at cell interfaces, i.e.,  $\bfv^{n+\hf}_{i+\hf,j}$ and $\bfw^{n+\hf}_{i,j+\hf}$, achieving second-order accuracy in time.

The mid-point values $\bfv^{n+\hf}_{i+\hf,j}$ and $\bfw^{n+\hf}_{i,j+\hf}$ are obtained by solving the GRP of the relaxation model \eqref{CNS2} at $(x_{i+\hf},y_j,t^n)$. This GRP is formulated as the IVP with piecewise smooth initial data defined on the two neighboring computational cells:
\begin{equation}
	\left.\begin{cases}
		\begin{aligned}
			&\frac{\pt\bfu}{\pt t}+\frac{\pt\bfv}{\pt x}+\frac{\pt\bfw}{\pt y}=0, \\
			&\frac{\pt\bfv}{\pt t}+a^2 \frac{\pt\bfu}{\pt x}
			=\frac{1}{\ep}\big[\bff_c(\bfu)-\bff_v(\bfu,\nabla\bfu)-\bfv\big],\\
			&\frac{\pt\bfw}{\pt t}+a^2 \frac{\pt\bfu}{\pt y}
			=\frac{1}{\ep}\big[\bfg_c(\bfu)-\bfg_v(\bfu,\nabla\bfu)-\bfw\big],\\
			&(\bfu,\bfv,\bfw)(\bfx,t^n)=\left.\begin{cases}
				(\bfu^n_{i,j},\bfv^n_{i,j},\bfw^n_{i,j})(\bfx), &\quad x<x_{i+\hf},\\
				(\bfu^n_{i+1,j},\bfv^n_{i+1,j},\bfw^n_{i+1,j})(\bfx),&\quad x>x_{i+\hf}.
			\end{cases}	
			\right.
		\end{aligned}
	\end{cases}\right.
	\label{FVGRP}
\end{equation}
In the above system, $a>0$ and $\ep>0$ are locally defined constant parameters which are determined based on the space data reconstruction at the cell interface $\Gamma_{i+\hf,j}=\{x_{i+\hf}\}\times[y_{j-\hf},y_{j+\hf}]$. Details will be shown in the next section. 

The piecewise smooth initial data of the conservative variable, i.e., $\bfu^n_{i,j}(\bfx)$ and $\bfu^n_{i+1,j}(\bfx)$, are obtained by the space data reconstruction procedure presented in \ref{datarecon}. The initial data of relaxation variables are given by the equilibrium state as
\begin{equation}\label{eq:equilibrium}
\bfv(\bfx,t^n)=\left.\begin{cases}
				\bff_c(\bfu^{n}_{i,j}(\bfx))-\bff_v(\bfu^{n}_{i,j}(\bfx),\nabla\bfu^{n}_{i,j}(\bfx)), &\quad x<x_{i+\hf},\\
				\bff_c(\bfu^{n}_{i+1,j}(\bfx))-\bff_v(\bfu^{n}_{i+1,j}(\bfx),\nabla\bfu^{n}_{i+1,j}(\bfx)),&\quad x>x_{i+\hf}.
			\end{cases}	
			\right.
\end{equation}

Solving the GRP \eqref{FVGRP} is the topic of the next section.
Here, we provide heuristic formulae to elucidate the essential idea of the GRP solver.
The mid-point value of $\bfv$ is approximated as
\begin{equation}\notag
\bfv^{n+\hf}_{i+\hf,j}=\bfv^{n,*}_{i+\hf,j}+\frac{\Delta t}{2}\bigg(\frac{\partial\bfv}{\partial t}\bigg)^{\text{impl},*}_{i+\hf,j},
\end{equation}
where $\bfv^{n,*}_{i+\hf,j}$ is the Riemann solution, and $\big(\frac{\partial\bfv}{\partial t}\big)^{\text{impl},*}_{i+\hf,j}$ stands for an estimation of the time derivative of $\bfv$. 
By performing the LW procedure to the governing PDEs of $\bfv$, the time derivative $\frac{\partial\bfv}{\partial t}$ takes the form
\begin{equation}\label{eq:LW-bfv}
	\frac{\partial\bfv}{\partial t}=-a^2 \frac{\pt\bfu}{\pt x}+\frac{1}{\ep}\big[\bff_c(\bfu)-\bff_v(\bfu,\nabla\bfu)-\bfv\big].
\end{equation}
Different from the approach for solving non-stiff hyperbolic PDEs, e.g., the Euler equations \cite{BLW2006,BL2007}, where $(\frac{\partial\bfv}{\partial t})^n$ is used, a more subtle implicit-explicit method is used here to cope with the stiffness of the governing PDEs of $\bfv$.
The first term on the right-hand side of \eqref{eq:LW-bfv} is the non-stiff one and the second term is the stiff one. The non-stiff and stiff terms are treated explicitly and implicitly, respectively, leading to
	\begin{equation}
		\bfv^{n+\hf}_{i+\hf,j}=\bfv_{i+\hf,j}^{n,*}
		+\frac{\Delta t}{2}
		\left[
		-a^2\bigg(\frac{\pt \bfu}{\pt x}\bigg)^{n,*}_{i+\hf,j}
		+\frac{1}{\epsilon}\bfH_{i+\hf,j}^{n+\hf}
		-\frac{1}{\epsilon}\bfv^{n+\hf}_{i+\hf,j}
		\right],
	\end{equation}
where 
	\begin{equation}
		\bfH=\bff_c(\bfu)-\bff_v(\bfu,\nabla\bfu). \label{Hdef}
	\end{equation}
	Consequently, $\bfv^{n+\hf}_{i+\hf,j}$ is expressed as
	\begin{equation*}
		\bfv^{n+\hf}_{i+\hf,j}=\frac{2\ep}{2\ep+\Delta t}\bfv_{i+\hf,j}^{n,*}
		-a^2\frac{\ep\Delta t}{2\ep+\Delta t}
		\bigg(\frac{\pt \bfu}{\pt x}\bigg)^{n,*}_{i+\hf,j}
		+\frac{\Delta t}{2\ep+\Delta t}\bfH_{i+\hf,j}^{n+\hf}.
	\end{equation*}
In the above formula, $\bfH_{i+\hf,j}^{n+\hf}$ is approximated by using the Crank-Nicolson method.
Therefore, the dependence of $\bfH$ on $\nabla\bfu$ leads to the necessity of solving linear systems of algebraic equations to update the conservative variable $\bfu$.

Due to the rotational invariance of the flux, the numerical flux $\bfw_{i,j+\hf}^{n+\hf}$ can be computed in the same manner.

\subsection{Set-up of the generalized Riemann problem}\label{sec:setup-grp}
The numerical flux $\bfv_{i+\hf,j}^{n+\hf}$ is obtained by solving the GRP at the cell interface $\Gamma_{i+\hf,j}$. Instead of directly solving the 2-D GRP \eqref{FVGRP}, we regard the transversal derivatives as source terms and formulate the quasi 1-D GRP \cite{LD2016}
\begin{equation}
	\left.\begin{cases}
		\begin{aligned}
			&\frac{\pt\bfu}{\pt t}+\frac{\pt\bfv}{\pt x}=-\frac{\pt\bfw}{\pt y}, \\
			&\frac{\pt\bfv}{\pt t}+a^2 \frac{\pt\bfu}{\pt x}
			=\frac{1}{\ep}\big[\bff_c(\bfu)-\bff_v(\bfu,\nabla\bfu)-\bfv\big],\\
			&(\bfu,\bfv)(x,y_j,t^n)=\left.\begin{cases}
				(\bm{\mathfrak{u}}_{L},\bm{\mathfrak{v}}_{L})(x,y_j), &\quad x<x_{i+\hf},\\
				(\bm{\mathfrak{u}}_{R},\bm{\mathfrak{v}}_{R})(x,y_j), &\quad x>x_{i+\hf}.
			\end{cases}	
			\right.
		\end{aligned}
	\end{cases}\right.
	\label{GRP}
\end{equation}
The initial data are piecewise linearly reconstructed and take the form
\begin{equation}
\begin{aligned}
    &\bm{\mathfrak{u}}_{L/R}(x,y)=\bfu_{L/R}+(\pt_x\bfu)_{L/R}\cdot(x-x_{i+\hf})+(\pt_y\bfu)_{L/R}\cdot(y-y_j),\\
&\bm{\mathfrak{v}}_{L/R}(x,y)=\bfv_{L/R}+(\pt_x\bfv)_{L/R}\cdot(x-x_{i+\hf}),
\end{aligned}
	\label{GRPini}
\end{equation}
where the notations $(\cdot)_L$ and $(\cdot)_R$ denote the limiting values on the left and right sides of the cell interface $\Gamma_{i+\hf,j}$, respectively.

For the relaxation variable $\bfv$, the constant vectors $\bfv_{L/R}$ and $(\pt_x\bfv)_{L/R}$ are defined as 
\begin{equation*}
	\begin{aligned}
		& \bfv_L=\bff_c(\bfu_L)-\bff_v(\bfu_L,(\nabla\bfu)_L),\quad
		\bfv_R=\bff_c(\bfu_R)-\bff_v(\bfu_R,(\nabla\bfu)_R),\\
		&(\pt_x\bfv)_L=(\pt_x\bff_c)_L-(\pt_x\bff_v)_L,\quad
		(\pt_x\bfv)_R=(\pt_x\bff_c)_R-(\pt_x\bff_v)_R,	
	\end{aligned}
\end{equation*}
where
\begin{equation}\notag
		(\pt_x\bff_c)_{L/R}=\frac{\pt \bff_c}{\pt\bfu}(\bfu_{L/R})(\pt_x\bfu)_{L/R},\quad
        (\pt_x\bff_v)_{L/R}
        =\left.\frac{\pt \bff_v(\bm{\mathfrak{u}}_{L/R},\nabla \bm{\mathfrak{u}}_{L/R})}{\pt x}\right|_{(x_{i+\hf},y_j)}.
\end{equation}

Since the initial data $\bm{\mathfrak{v}}_L(x,y)$ and $\bm{\mathfrak{v}}_R(x,y)$ are solely determined by $\bm{\mathfrak{u}}_L(x,y)$ and $\bm{\mathfrak{u}}_R(x,y)$, there is no need for the finite volume scheme \eqref{FV} to evolve $\bfv$. The reconstruction procedure of $\bfv$ is also avoided. The same conclusion holds for $\bfw$ as well. Therefore, the relaxation model \eqref{CNS2} is only used to compute numerical fluxes for the conservative variable $\bfu$. An additional advantage of the present numerical scheme is that the parameters $\ep$ and $a$ 
 are not necessarily to be universal either in space or in time in numerical computations. 

\section{Generalized Riemann problem solver for the relaxation model} \label{GRPFS}
This section focuses on solving the GRP to compute $\bfv_{i+\hf,j}^{n+\hf}$ and $\bfw_{i,j+\hf}^{n+\hf}$, which are essential components of the finite volume scheme \eqref{FV}. 
With a sufficiently small relaxation parameter $\ep$, the numerical procedure developed in this section can be regarded as an approximate GRP solver for the compressible Navier-Stokes equations, which is termed RFS.

\subsection{Associated Riemann solver}
The first step of solving the GRP is to compute the Riemann solution at the cell interface.
For the IVP \eqref{GRP} of the linear hyperbolic equation system, we employ the treatment of diagonalization and tracking the characteristics. Without any loss of generality, we shift the coordinates such that $(x_{i+\hf},y_j,t^n)=(0,0,0)$.

Rewrite the governing equations in \eqref{GRP} as
\begin{equation}
	\left.\begin{cases}
		\begin{aligned}
			&\frac{\pt(\bfv+a\bfu)}{\pt t} 
			+ a\frac{\pt(\bfv+a\bfu)}{\pt x}
			=-a\frac{\pt\bfw}{\pt y}+\frac{1}{\epsilon}(\bfH-\bfv), \\
			&\frac{\pt(\bfv-a\bfu)}{\pt t} 
			- a\frac{\pt(\bfv-a\bfu)}{\pt x}
			=a\frac{\pt\bfw}{\pt y}+\frac{1}{\epsilon}(\bfH-\bfv).
		\end{aligned}
	\end{cases}\right.
	\label{charaline3}
\end{equation}
Integrate the equations $\eqref{charaline3}$ along the two families of characteristic curves ending at $(\bfx_0,t)=(0,0,\theta\Delta t)$ for $\theta\in[0,1]$, respectively. This leads to the result
\begin{equation}
	\begin{aligned}
		&\bfv(0,0,\theta\Delta t) \pm a \bfu(0,0,\theta\Delta t)\\
		=&\bm{\mathfrak{v}}_{J_\pm}(\mp a\theta\Delta t,0) \pm a \bm{\mathfrak{u}}_{J_\pm}(\mp a\theta\Delta t,0)
		\mp a\int_0^{\theta\Delta t}\frac{\pt\bfw}{\pt y}(\pm a(s-\theta\Delta t),0,s)ds\\
		&+\frac{1}{\ep}\int_0^{\theta\Delta t} \bfH(\pm a(s-\theta\Delta t),0,s)
		-\bfv(\pm a(s-\theta\Delta t),0,s)ds,
	\end{aligned}
	\label{charaline4}
\end{equation}
where $J_+=L$ and $J_-=R$. Then, we approximate the integrals in \eqref{charaline4} by the forward and backward Euler methods for non-stiff and stiff source terms, respectively. As a result,
	\begin{equation}
		\begin{aligned}
			\bfv^\theta_{\bfx_0}\pm a\bfu^\theta_{\bfx_0}
			=&\bm{\mathfrak{v}}_{J_\pm}(\mp a\theta\Delta t,0) \pm a\bm{\mathfrak{u}}_{J_\pm}(\mp a\theta\Delta t,0)\\
			&-a\theta\Delta t(\pt_y \bfw)_{J_\pm}+\frac{\theta\Delta t}{\ep}(\bfH^\theta_{\bfx_0}-\bfv^\theta_{\bfx_0}),
		\end{aligned}
		\label{GRPsol1}
	\end{equation}
where $\bfv^\theta_{\bfx_0}$, $\bfu^\theta_{\bfx_0}$ and $\bfH^\theta_{\bfx_0}$ denote the approximations to $\bfv$, $\bfu$ and $\bfH$ at $(0,0,\theta \Delta t)$ locating on the cell interface, and $(\pt_y\bfw)_{J_\pm}$ is the constant vector approximating $\frac{\pt \bfw}{\pt y}$ at the starting point $(0\mp a\theta\Delta t,0,0)$.
Similarly to \eqref{GRPini}, $(\pt_y\bfw)_L$ and $(\pt_y\bfw)_R$ are defined as
\begin{equation}
		(\pt_y\bfw)_{L/R}=(\pt_y\bfg_c)_{L/R}-(\pt_y\bfg_v)_{L/R},
\end{equation}
where
\begin{equation}\notag
		(\pt_y\bfg_c)_{L/R}=\frac{\pt \bfg_c}{\pt\bfu}(\bfu_{L/R})(\pt_y\bfu)_{L/R},\quad
        (\pt_y\bfg_v)_{L/R}
        =\left.\frac{\pt \bfg_v(\bm{\mathfrak{u}}_{L/R},\nabla \bm{\mathfrak{u}}_{L/R})}{\pt y}\right|_{(0,0)}.
\end{equation}

It follows from \eqref{GRPsol1} that, for $\theta\in[0,1]$,
\begin{equation}
	\begin{aligned}
		\bfu^{\theta}_{\bfx_0}
		=&\frac{1}{2}\big[\bm{\mathfrak{u}}_L(-a\theta\Delta t,0)+\bm{\mathfrak{u}}_R(a\theta\Delta t,0)\big]
		-\frac{1}{2a}\big[\bm{\mathfrak{v}}_R(a\theta\Delta t,0)-\bm{\mathfrak{v}}_L(-a\theta\Delta t,0)\big]\\
		&-\frac{1}{2}\theta\Delta t\big[(\pt_y\bfw)_R-(\pt_y\bfw)_L\big],\\
		\bfv^{\theta}_{\bfx_0}
		=&\bigg\{\frac{1}{2}\big[\bm{\mathfrak{v}}_L(-a\theta\Delta t,0)+\bm{\mathfrak{v}}_R(a\theta\Delta t,0)\big]
		-\frac{a}{2}\big[\bm{\mathfrak{u}}_R(a\theta\Delta t,0)-\bm{\mathfrak{u}}_L(-a\theta\Delta t,0)\big]\\
		&-\frac{a}{2}\theta\Delta t\big[(\pt_y\bfw)_R-(\pt_y\bfw)_L\big]
		+\frac{\theta\Delta t}{\ep} \bfH^{\theta}_{\bfx_0}\bigg\}
		\bigg/\bigg(1+\frac{\theta\Delta t}{\ep}\bigg).
	\end{aligned}
	\label{GRPsol2}
\end{equation}
The Riemann solution $\big(\bfu^{n,*}_{i+\hf,j}, \bfv^{n,*}_{i+\hf,j}\big)$ is consequently obtained by substituting $\theta=0$ into \eqref{GRPsol2} and it holds that,
\begin{subequations}\label{RiemannSol}
\begin{equation}
	\bfu_{i+\hf,j}^{n,*}=\frac{1}{2}\big(\bfu_L+\bfu_R\big) 
	-\frac{1}{2a}\big(\bfv_R-\bfv_L\big),   \label{GRPsol3}
\end{equation}
\begin{equation}
	\bfv_{i+\hf,j}^{n,*}=\frac{1}{2}\big(\bfv_L+\bfv_R\big)
	-\frac{a}{2}\big(\bfu_R-\bfu_L\big).  \label{GRPsol4}
\end{equation}
\end{subequations} 

Following the terminology of \cite{BF2003,BLW2006}, we term the numerical procedure to obtain \eqref{RiemannSol} the associated Riemann solver.

\subsection{Generalized Riemann problem solver} \label{mid-point}

With the formula \eqref{GRPsol2} at hand, we are in a position to compute the mid-point cell interface value $\bfv^{n+\hf}_{i+\hf,j}$ used in \eqref{FV}.

By setting $\theta=\frac{1}{2}$ in \eqref{GRPsol2}, we obtain, after using the initial condition \eqref{GRPini}, that
\begin{equation}	\bfv^{n+\hf}_{i+\hf,j}=\frac{2\ep}{2\ep+\Delta t}\bfv_{i+\hf,j}^{n,*}
		-a^2\frac{\ep\Delta t}{2\ep+\Delta t}
		\bigg(\frac{\pt \bfu}{\pt x}\bigg)^{n,*}_{i+\hf,j}
		+\frac{\Delta t}{2\ep+\Delta t}\bfH_{i+\hf,j}^{n+\hf},
		\label{Taylor}
\end{equation}
where the first term on the right-hand side of the above formula is already obtained by \eqref{GRPsol4}, and the second term is expressed as
\begin{equation}
\begin{aligned}
    \bigg(\frac{\pt \bfu}{\pt x}\bigg)^{n,*}_{i+\hf,j}
	=&\frac{1}{2}\big[(\pt_x\bfu)_L+(\pt_x\bfu)_R\big]-\frac{1}{2a}\big[(\pt_x\bfv)_R-(\pt_x\bfv)_L \big]\\
	&-\frac{1}{2a}\big[(\pt_y\bfw)_R-(\pt_y\bfw)_L \big].   \label{GRPsol5}
\end{aligned}
\end{equation}
To achieve second-order accuracy in time, the last term on the right-hand side of \eqref{Taylor}, $\bfH_{i+\hf,j}^{n+\hf}$, is computed in the spirit of the Crank-Nicolson method, i.e., 
	\begin{equation} \label{eq:H-Crank-Nicolson}
		\bfH_{i+\hf,j}^{n+\hf}
		=\frac{1}{2}\big(\bfH_{i+\hf,j}^{n}+\bfH_{i+\hf,j}^{n+1,-}\big),
	\end{equation}	
where
\begin{equation}\notag
	\bfH_{i+\hf,j}^{n}
	=(\bff_c)_{i+\hf,j}^{n}-(\bff_v)_{i+\hf,j}^{n},\quad
	\bfH_{i+\hf,j}^{n+1,-}
	=(\bff_c)_{i+\hf,j}^{n+1,-}-(\bff_v)_{i+\hf,j}^{n+1,-}.
\end{equation}
At $t=t^n$, the convective part $(\bff_c)_{i+\hf,j}^{n}$ is directly obtained by using the Riemann solution \eqref{GRPsol3},
\begin{equation}\notag
(\bff_c)_{i+\hf,j}^{n}=\bff_c(\bfu_{i+\hf,j}^{n,*}).
\end{equation}
To compute the viscous part $(\bff_v)_{i+\hf,j}^{n}$, we express the viscous fluxes by regarding the primitive variables $\bfQ(\bfu)=(\rho,u,v,T)^\top$ as independent variables, i.e.,
\begin{equation*}
	\bff_v(\bfu,\nabla\bfu)=\bftf_v(\bfQ,\nabla\bfQ),\quad
	\bfg_v(\bfu,\nabla\bfu)=\bftg_v(\bfQ,\nabla\bfQ).
\end{equation*}
So, $(\bff_v)_{i+\hf,j}^{n}$ is expressed as
\begin{equation}\label{eq:viscous-flux-n}
(\bff_v)_{i+\hf,j}^{n}
	=\bftf_v\big(\bfQ_{i+\hf,j}^{n},(\nabla\bfQ)_{i+\hf,j}^{n}\big)
\end{equation}
where
\begin{equation}
	\begin{aligned}
&\bfQ_{i+\hf,j}^{n}=\bfQ(\bfu_{i+\hf,j}^{n,*}),
    \\
	&\bigg(\frac{\pt\bfQ}{\pt x}\bigg)_{i+\hf,j}^{n}= \frac{(\pt_x\bfQ)^{n}_{L}  + (\pt_x\bfQ)^{n}_{R}}{2}+\frac{(\bfQ^{n}_{R}-\bfQ^{n}_{L})}{\Delta x_{i+\hf}},\\
    \quad
    &\bigg(\frac{\pt\bfQ}{\pt y}\bigg)_{i+\hf,j}^{n}
	=\frac{(\pt_y\bfQ)^{n}_{L}+(\pt_y\bfQ)^{n}_{R}}{2},
	\end{aligned}
	\label{gradient1-0}
\end{equation}
and $\Delta x_{i+\hf}=\frac{1}{2}(\Delta x_i+\Delta x_{i+1})$.
Here, $(\pt_{x}\bfQ)^{n}_{L/R}$ and $(\pt_{y}\bfQ)^{n}_{L/R}$ come from the space date reconstruction and details are put in \ref{datarecon}.
The approximation $(\frac{\pt\bfQ}{\pt x})_{i+\hf,j}^{n}$ is obtained in the spirit of the direct discontinuous Galerkin (DDG) method \cite{LY2009}.

At $t=t^{n+1}$, the convective part $(\bff_c)_{i+\hf,j}^{n+1,-}$ is
\begin{equation*}
(\bff_c)_{i+\hf,j}^{n+1,-}=\bff_c(\bfu_{i+\hf,j}^{n+1,-}),
\end{equation*}
where $\bfu_{i+\hf,j}^{n+1,-}$ is obtained by setting $\theta=1$ in  \eqref{GRPsol2},
\begin{equation}
	\begin{aligned}
\bfu_{i+\hf,j}^{n+1,-}
=&\bfu_{i+\hf,j}^{n,*}
    +\frac{\Delta t}{2}\big[(\pt_x\bfv)_R+(\pt_x\bfv)_L \big] \\
    &-\frac{a\Delta t}{2}\big[(\pt_x\bfu)_R-(\pt_x\bfu)_L\big]
	+\frac{\Delta t}{2}\big[(\pt_y\bfw)_R+(\pt_y\bfw)_L \big].
	\end{aligned}
	\label{GRPsol6}
\end{equation}
Similarly to \eqref{eq:viscous-flux-n}, the viscous part at $t^{n+1}$ is
\begin{equation}
	(\bff_v)_{i+\hf,j}^{n+1,-}
	=\bftf_v\big(\bfQ_{i+\hf,j}^{n+1,-},(\nabla\bfQ)_{i+\hf,j}^{n+1,-}\big),
	\label{visterm}
\end{equation}
where
\begin{equation}
	\begin{aligned}
&\bfQ_{i+\hf,j}^{n+1,-}=\bfQ(\bfu_{i+\hf,j}^{n+1,-}), \\
& \bigg(\frac{\pt\bfQ}{\pt x}\bigg)_{i+\hf,j}^{n+1,-}
= \frac{(\pt_x\bfQ)^{n+1,-}_{L}
+(\pt_x\bfQ)^{n+1,-}_{R}}{2}
+\frac{\bfQ^{n+1,-}_{R}-\bfQ^{n+1,-}_{L}}{\Delta x_{i+\hf}},\\
&\bigg(\frac{\pt\bfQ}{\pt y}\bigg)_{i+\hf,j}^{n+1,-}
=\frac{(\pt_y\bfQ)^{n+1,-}_{L}+(\pt_y\bfQ)^{n+1,-}_{R}}{2}.
	\end{aligned}
	\label{gradient1-1}
\end{equation}
The values of $(\pt_x\bfQ)^{n+1,-}_{,L/R}$ and $(\pt_y\bfQ)^{n+1,-}_{,L/R}$ are obtained by
\begin{equation}\notag
	\begin{aligned}
&(\pt_x\bfQ)^{n+1,-}_{L}=\bfbsig_{i,j}^{n+1},
\quad
(\pt_x\bfQ)^{n+1,-}_{R}=\bfbsig_{i+1,j}^{n+1},\\
&(\pt_y\bfQ)^{n+1,-}_{L}=\bfbvpi_{i,j}^{n+1},
\quad
(\pt_y\bfQ)^{n+1,-}_{R}=\bfbvpi_{i+1,j}^{n+1},
	\end{aligned}
\end{equation}
where $(\bfbsig_{i,j}^{n+1},\bfbvpi_{i,j}^{n+1})$ stands for the gradient estimation of $\bfQ$ over $\Omega_{i,j}$ at $t=t^{n+1}$ by the GRP methodology \cite{BLW2006}. 
Using $\bfQ^{n+1,-}_{i\pm\hf,j}=\bfQ(\bfu^{n+1,-}_{i\pm\hf,j})$ and $\bfQ^{n+1}_{i,j\pm\hf}=\bfQ(\bfu^{n+1,-}_{i,j\pm\hf})$, they are computed as
\begin{equation}
	\bfbsig_{i,j}^{n+1}=\frac{1}{\Delta x_i}
	(\bfQ_{i+\hf,j}^{n+1,-}-\bfQ_{i-\hf,j}^{n+1,-}),\quad
	\bfbvpi_{i,j}^{n+1}=\frac{1}{\Delta y_j}
	(\bfQ_{i,j+\hf}^{n+1,-}-\bfQ_{i,j-\hf}^{n+1,-}).
	\label{slope}
\end{equation} 
The values of $\bfQ^{n+1,-}_{R}$ and $\bfQ^{n+1,-}_{L}$ are approximated with linear interpolation as
\begin{equation}
\bfQ^{n+1,-}_{L}
=\bftQ_{i,j}^{n+1}
	+\frac{\Delta x_{i}}{2}\bfbsig_{i,j}^{n+1},\quad
\bfQ^{n+1,-}_{R}
=\bftQ_{i+1,j}^{n+1}
	-\frac{\Delta x_{i+1}}{2}\bfbsig_{i+1,j}^{n+1}, \label{implicit0}
\end{equation} 
where $\bftQ_{i,j}^{n+1}$ depends on conservative variables through
\begin{equation}
\bftQ_{i,j}^{n+1}=\bfQ(\bfbu_{i,j}^{n+1}).
\label{implicitness}
\end{equation} 

Substituting \eqref{Taylor}-\eqref{implicitness} into the numerical scheme \eqref{FV} leads to the presence of $\bfbu_{i,j}^{n+1}$ on both sides of the equation, resulting in systems of algebraic equations to solve. In addition, the systems are discovered to be linear.

\begin{rem}
The implicit property of the RFS is due to the use of \eqref{gradient1-1} in the computation of $\bfH_{i+\hf,j}^{n+\hf}$. With the interface value $\bfQ_{i+\hf,j}^{n+1,-}$ already provided by the GRP solver through \eqref{GRPsol6}, the implicit treatment only involves spatial derivatives used in \eqref{visterm}. Then, the DDG style approximation of $(\nabla\bfQ)_{i+\hf,j}^{n+1,-}$ in \eqref{gradient1-1} results in linear systems of algebraic equations. Such linear systems will be formulated in the next subsection. 
\end{rem}

\subsection{Solving the linear system of algebraic equations}

This subsection formulates the linear systems of algebraic equations to solve for the update of the cell average $\bar{\bfu}_{i,j}^{n+1}$ in the finite volume scheme \eqref{FV}.

We first note that the cell average of the density $\bar\rho^{n+1}_{i,j}$ can be updated through \eqref{FV} explicitly. 
After substituting $\bfv^{n+\hf}_{i+\hf,j}$ and $\bfw^{n+\hf}_{i,j+\hf}$ into the numerical scheme \eqref{FV}, it is derived from the momentum equation in the $x$-direction that
\begin{equation}\label{eq:impl-evol-u}
\begin{aligned}
&\displaystyle\frac{\bar\rho_{i,j}^{n+1}\tilde{u}_{i,j}^{n+1}-\bar\rho_{i,j}^{n}\tilde{u}_{i,j}^{n}}{\Delta t}
    -\sum_{s=1}^{6}R_{i,j,s}\\
=&\frac{2}{3}\frac{\mu(T^{n+1,-}_{i+\hf,j})\omega^n_{i,j+\hf}\frac{\tilde{u}^{n+1}_{i+1,j}-\tilde{u}^{n+1}_{i,j}}{\Delta x_{i+\hf}}-
\mu(T^{n+1,-}_{i-\hf,j})\omega^n_{i-\hf,j}\frac{\tilde{u}^{n+1}_{i,j}-\tilde{u}^{n+1}_{i-1,j}}{\Delta x_{i-\hf}}}{\Delta x_i}\\
&+\frac{1}{2}\frac{\mu(T^{n+1,-}_{i,j+\hf})\omega^n_{i,j+\hf}\frac{\tilde{u}^{n+1}_{i,j+1}-\tilde{u}^{n+1}_{i,j}}{\Delta y_{j+\hf}}
-\mu(T^{n+1,-}_{i,j-\hf})\omega^n_{i,j-\hf}\frac{\tilde{u}^{n+1}_{i,j}-\tilde{u}^{n+1}_{i,j-1}}{\Delta y_{j-\hf}}}{\Delta y_j},
\end{aligned}
\end{equation}
where
\begin{equation*}
		\omega^n_{i+\hf,j} = \frac{\Delta t}{2\ep^n_{i+\hf,j}+\Delta t },
		\quad \omega^n_{i,j+\hf}  = \frac{\Delta t}{2\ep^n_{i,j+\hf}+\Delta t},
\end{equation*}
and $R^{u}_{i,j,s}$ are the explicit parts of numerical fluxes (see \ref{terms} for details). As stated in \Cref{sec:setup-grp}, 
the parameters, $\ep$ and $a$,
are locally defined here.

By grouping evolutions \eqref{eq:impl-evol-u} for all $i\in\{0,1,\dots,I\}$ and $j\in\{0,1,\dots,J\}$, we obtain a linear equation system regarding $\tilde{u}_{i,j}^{n+1}$. Linear equation systems regarding $\tilde{v}_{i,j}^{n+1}$ and $\tilde{T}_{i,j}^{n+1}$ can be established in the same manner.
Denote linear systems as
\begin{equation} \label{linear-sys}
	\fM_{u}\vu=\vb_u,\quad 	\fM_{v}\vv=\vb_v, \quad	\fM_{T}\vT=\vb_T.
\end{equation}
The unknown vector for the velocity component $u$ is
\begin{equation*}
	\vu=\begin{pmatrix}
		\tilde{u}^{n+1}_{1,1},\tilde{u}^{n+1}_{1,2},\dots,\tilde{u}^{n+1}_{1,J},
		\tilde{u}^{n+1}_{2,1}, \tilde{u}^{n+1}_{2,2},\dots,\tilde{u}^{n+1}_{2,J},\dots,
		\tilde{u}^{n+1}_{I,1}, \tilde{u}^{n+1}_{I,2},\dots,\tilde{u}^{n+1}_{I,J}
	\end{pmatrix}^\top.
\end{equation*} 
For $k= (i-1)\times I+j$, the $(k,l)$-th entry of $\fM_{u}$ is defined as
\begin{equation} \label{Jacobian-u}
	(\fM_u)_{k,l}=\begin{cases}
		\frac{4}{3}\mu(T^{n+1,-}_{i-\hf,j})\frac{\Delta t}{2\Delta x_i \Delta x_{i-\hf}}\omega^n_{i-\hf,j},
		\quad\quad &\text{if}\,\, l=k-I,\\[2mm]
		\frac{4}{3}\mu(T^{n+1,-}_{i+\hf,j})\frac{\Delta t}{2\Delta x_i \Delta x_{i+\hf}}\omega^n_{i+\hf,j},
		\quad\quad &\text{if}\,\, l=k+I,\\[2mm]  
		\mu(T^{n+1,-}_{i,j-\hf})\frac{\Delta t}{2\Delta y_j \Delta y_{j-\hf}}\omega^n_{i,j-\hf},
		\quad\quad &\text{if}\,\, l=k-1,\\[2mm]
		\mu(T^{n+1,-}_{i,j+\hf})\frac{\Delta t}{2\Delta y_j \Delta y_{j+\hf}}\omega^n_{i,j+\hf},
		\quad\quad &\text{if}\,\, l=k+1,\\[2mm] 
		\bar\rho^{n+1}_{i,j}
		+(\fM_u)^k_{k-I}+(\fM_u)^k_{k+I}\\
		+(\fM_u)^k_{k-1}+(\fM_u)^k_{k+1},     &\text{if}\,\, l=k,\\[2mm]
		0, \quad\quad\quad\quad          &\text{otherwise}.
	\end{cases}	         
\end{equation} 

For an interior cell in the computational domain, the $k$-th entry of the right-hand side vector $\vb_u$ is
\begin{equation*}
(\vb_u)_k=\bar\rho_{i,j}^{n}\tilde{u}_{i,j}^{n}+\Delta t\sum_{s=1}^{6}R^{u}_{i,j,s},
\end{equation*}
where $(i,j)$ is the cell index such that $(i-1)\times I+j=k$.
Entries involving boundary cells in $\fM_u$ and $\vb_u$ are treated by to the conventional way of extrapolation in ghost cells and we do not elaborate here for the simplicity of the manuscript. The coefficient matrices and right-hand side vectors of the linear equation systems for $v$ and $T$ are put in \ref{terms}.

It is seen from \eqref{Jacobian-u} and \eqref{Jacobian-T} that the coefficient matrices $\fM_{u}$, $\fM_v$ and $\fM_T$ are sparse and diagonally-dominant. In the present work, the linear systems \eqref{linear-sys} are solved by the classical Jacobian iterative method.

\begin{rem}
The time stepping formula \eqref{eq:impl-evol-u} resembles the discretization to a diffusion equation of $u$. Solving the linear systems \eqref{linear-sys} acts just like the implicit finite difference discretization utilized in RK type schemes. The severe time step restriction imposed by the stability condition for diffusion terms is consequently eliminated.
This is also demonstrated numerically by a series of benchmark tests in \Cref{Examples}.  
\end{rem}

\subsection{Parameters}
This section presents the parameters used in the current work. All these parameters are locally determined based on the reconstructed initial data around cell interfaces.

The parameter $a$ in \eqref{GRP} is given by
\begin{equation}\label{eq:art-vis}
    a=\max(|u_L|+c_L, |u_R|+c_R),
\end{equation}
where $c=\sqrt{\gamma T}$ is the sound speed.

The relaxation parameter $\ep$ is set to be
\begin{equation}\label{eq:relaxation-param}
	\ep=\ep_0+C_{num}\frac{|p_R-p_L|}{p_R+p_L}\Delta t,
\end{equation}
where $\ep_0=10^{-9}$ and $C_{num}$ is a positive constant corresponding to the numerical viscosity. In the present paper, $C_{num}=1$ for all viscous flow problems, whlie $C_{num}=5$ for all inviscid flow problems due to the lack of physical viscosity.

The choice of $\ep$ is based on a similar consideration to the GKS method \cite{Xu2001}. In the shock wave region, the artificial dissipation brought by $\ep$ is at order $O(\Delta t)$ to improve the robustness. In the smooth flow region, the artificial dissipation brought by $\ep$ becomes very small such that the expected accuracy of the numerical method can be achieved.

\section{Summary of the algorithm}\label{sec:algorithm}

In this section, we summarize the numerical scheme developed so far as \Cref{alg:scheme}, including the finite volume scheme as well as the relaxation flux solver.

\begin{algorithm}[!htbp]
  \caption{Finite volume scheme using RFS} \label{alg:scheme}
  \begin{algorithmic}[1]
    \State Perform the reconstruction procedure to obtain $\bfu_{i,j}^n(\bfx)$ in each cell $\Omega_{i,j}$ for $i=1,\dots,I$ and $j=1,\dots,J$.
    \State At the cell interface $\Gamma_{i+\hf,j}$, set the values of the parameters $a$ and $\ep$ in the relaxation model according to \eqref{eq:art-vis} and \eqref{eq:relaxation-param}, respectively. Construct the initial data for $\bfv$ accoroding to \eqref{GRPini}.
    \State Apply the associated Riemann solver to compute $(\bfu_{i+\hf,j}^{n,*}, \bfv_{i+\hf,j}^{n,*})$ by \eqref{RiemannSol}. Compute $\bfu_{i+\hf,j}^{n+1,-}$ by \eqref{GRPsol6}.
    \State At the cell interface $\Gamma_{i,j+\hf}$, set $a$ and $\ep$, construct the initial data for $\bfw$, and compute $(\bfu_{i,j+\hf}^{n,*}, \bfw_{i,j+\hf}^{n,*})$ and $\bfu_{i,j+\hf}^{n+1,-}$ by the similar procedure.
    \State Estimate $(\bar{\bm{\sigma}}_{i,j}^{n+1},\bar{\bm{\varpi}}_{i,j}^{n+1})$, the slopes of $\bfQ$, in $\Omega_{i,j}$ by \eqref{slope}.
    \State Update $\bar{\rho}^{n+1}_{i,j}$ explicitly by \eqref{FV} and \eqref{Taylor}.
    \State For $u$, compute $\{R^u_{i,j,s}:s=1,\dots,6\}$ by \eqref{eq:def-R1u}-\eqref{eq:def-R456u}. 
    For $v$, compute $\{R^v_{i,j,s}:s=1,\dots,6\}$ similarly.
     For $T$, compute $\{R^T_{i,j,s}:s=1,\dots,6\}$ by \eqref{eq:def-R123T} and \eqref{eq:def-R456T}. 
    \State Solving linear systems \eqref{linear-sys} to update the rest of $\bar\bfu_{i,j}^{n+1}$.
  \end{algorithmic}
\end{algorithm}

In simulations for inviscid compressible flows governed by the Euler equations, the dynamic viscosity coefficient $\mu=0$, and the viscous fluxes $\bff_v(\bfu,\nabla \bfu)$ and $\bfg_v(\bfu,\nabla \bfu)$ vanish. Thus, $\bar{\bfu}^{n+1}_{i,j}$ is updated explicitly and there is no need to solve linear systems. Consequently, steps 7 and 8 in \Cref{alg:scheme} are omitted in computations for inviscid flows.  

\section{Numerical experiments} \label{Examples}
Numerical results of several benchmark problems, involving both inviscid and viscous compressible flows, are presented in this section to illustrate the performance of the current scheme.
In all the benchmarks, we apply the conventional CFL condition used in computations of the Euler equations, which takes the form
\begin{equation}\notag
	\Delta t=\lambda_{CFL}\frac{\min_{1\leq i\leq I, 1\leq j\leq J }(h_{i,j})}{\max_{1\leq i\leq I, 1\leq j\leq J }(|u_{i,j}|+|v_{i,j}|+c_{i,j})}.
\end{equation}
Here, $\lambda_{CFL}$ is the CFL number and $h_{i,j}$ is the cell size of the cell $\Omega_{i,j}$. 
We set $\lambda_{CFL}=0.6$ for 1-D cases and $\lambda_{CFL}=0.4$ for 2-D cases. The Prandtl number is $\Pr=0.72$ and the specific heat ratio is set to be $\gamma=1.4$. The data reconstruction procedure as shown in \ref{datarecon} is adopted, and the limiter parameter is taken as $\alpha=1.30$. To solve linear systems of algebraic equations, the Jacobian iteration method is employed. The boundary conditions are implemented by the extrapolation technique.

\subsection{Inviscid flows}

\vspace{2mm}
\subsubsection{Advection of a density perturbation}
The first test case is the advection of a density perturbation in 1-D
case. This case is used to validate the accuracy of the current scheme for inviscid smooth flows. The flow field is initialized by the analytical solution,
\begin{equation}
	\rho(x,t)=1+0.2\sin(\pi(x-t)),\quad u(x,t)=1,\quad p(x,t)=1, \notag
\end{equation}
with $t=0$. The computational domain is $\Omega =[0, 2]$ with periodic boundary conditions. A uniform mesh with $N$ cells is used. 

The $L^1$- and $L^\infty$-errors and convergence orders for the density $\rho$ at $t=2$ are presented in \Cref{Tab:Adden}. It is observed that the expected order of accuracy is achieved. 

\begin{table}[!htbp]
	\centering
	\caption{\footnotesize Accuracy test for the advection of a density perturbation: the $L^1$- and $L^\infty$-errors and convergence orders for the density $\rho$ at $t=2$.}
	\setlength{\tabcolsep}{2.5mm}
	\begin{tabular}{c|cc|cc}
		\toprule
		N &$L^1$ error&  convergence order  
		&$L^\infty$ error  & convergence order
		\\
		\hline
		
		20   &3.922e-03 &    &7.247e-03   &    
		\\
		40   &9.469e-04 &2.050   &1.667e-03   &2.120 
		\\
		80 &2.344e-04  &2.014   &3.984e-04   &2.065   
		\\
		160  &5.851e-05 &2.002   &9.826e-05   &2.020     
		\\
		320 &1.463e-05 &2.000    &2.450e-05   &2.004     
		\\
		640  &3.657e-06 &2.000   &6.119e-06   &2.002     
		\\
		1280  &9.143e-07 &2.000    &1.525e-06   &2.004    
		\\
		\bottomrule
	\end{tabular}   
    \label{Tab:Adden}  
\end{table}

\vspace{2mm}

\subsubsection{Isentropic vortex propagation problem}
The next accuracy test is the isentropic vortex propagation problem \cite{Shu1998}. The mean flow is $(\rho_0, u_0, v_0, p_0) = (1, 1, 1, 1)$, and an isotropic vortex is added to the mean flow, i.e., with perturbations added to the velocity $(u,v)$ and the temperature $T = \frac{p}{\rho}$, and no perturbation added to the entropy $S = \frac{p}{\rho^{\gamma}}$. The perturbation is given by
\begin{equation*}
	\begin{aligned}
		&(\delta u, \delta v)=\frac{\xi}{2\pi}e^{\frac{(1-r^2)}{2}}(y,x),\quad
		\delta T =-\frac{(\gamma-1)\xi^2}{8\gamma\pi^2}e^{1-r^2},
		\quad \delta S=0,
	\end{aligned}
\end{equation*}
where $r^2=x^2+y^2$ and $\xi=5$. The computational domain is $[-10,10]\times[-10,10]$, and the periodic boundary condition is imposed at all boundaries. The exact solution is the perturbation propagating with the background velocity.

The $L^1$- and $L^\infty$-errors and convergence orders for the density $\rho$ at $t=20$ with $N \times N$ uniform mesh cells are presented in \Cref{Tab:IVP}, which validate that the expected accuracy is achieved for the 2-D inviscid flow.

\begin{table}[!htbp]
	\centering
	\caption{\footnotesize Accuracy test for the isotropic vortex propagation problem: the $L^1$- and $L^\infty$-errors and convergence orders for the density $\rho$ at $t=20$.}   \label{Tab:IVP}  
	\setlength{\tabcolsep}{2.5mm}
	\begin{tabular}{c|cc|cc}
		\toprule
		N &$L^1$ error&   order 
		&$L^\infty$ error  &  order
		\\
		\hline
		$40$   &3.822e-03 &   &2.747e-01  &   
		\\
		$80$ &1.116e-03 &1.776   &1.082e-01  &1.344    
		\\
		$160$  &2.205e-04 &2.339   &2.184e-02   &2.308    
		\\
		$320$  &4.865e-05 &2.180   &4.484e-03   &2.284    
		\\
		$640$  &1.186e-05 &2.035   &1.003e-03   &2.160    
		\\
		$1280$   &2.959e-06 &2.003  &2.376e-04  &2.077
		\\
		\bottomrule
	\end{tabular}   
\end{table}  

\vspace{2mm}

\subsubsection{2-D Riemann problems}
We proceed to verify the performance of the current scheme in capturing multi-dimensional wave structures by considering two Riemann problems.
The selected two 2-D Riemann problems involve the interaction of planar rarefaction waves \cite{PXL2017}
and the interaction of vortex-sheets with different signs \cite{PXL2017}.
The computational domain is $[0,1]\times[0,1]$, and a uniform mesh with $1200\times 1200$ cells is used. Along the boundaries, the flow field is extended into ghost cells with interior values.

\begin{figure}[!htb]
	\centering
	\subfigure{	\includegraphics[width=0.45\linewidth]{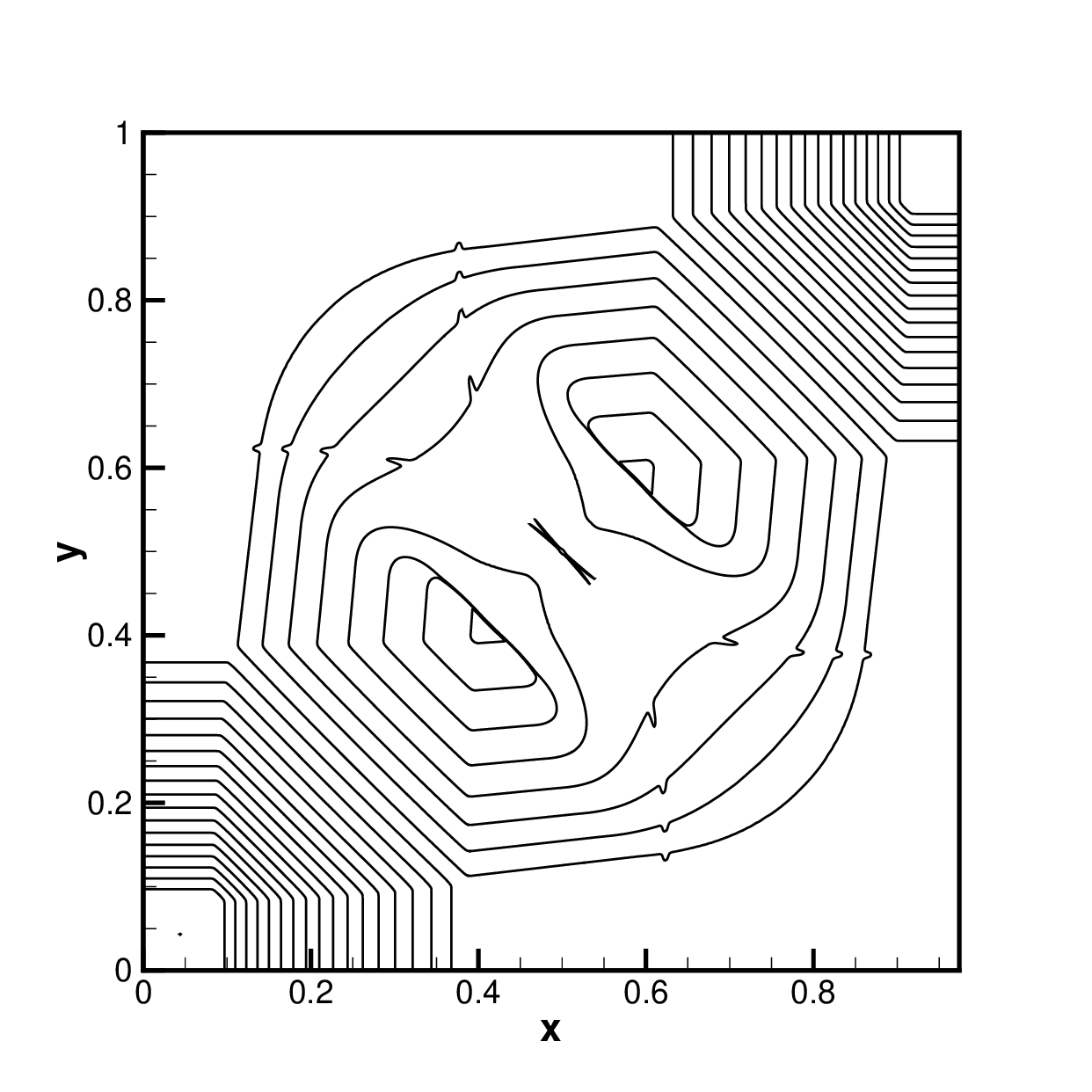}}
	\hspace{0mm}
	\subfigure{	\includegraphics[width=0.45\linewidth]{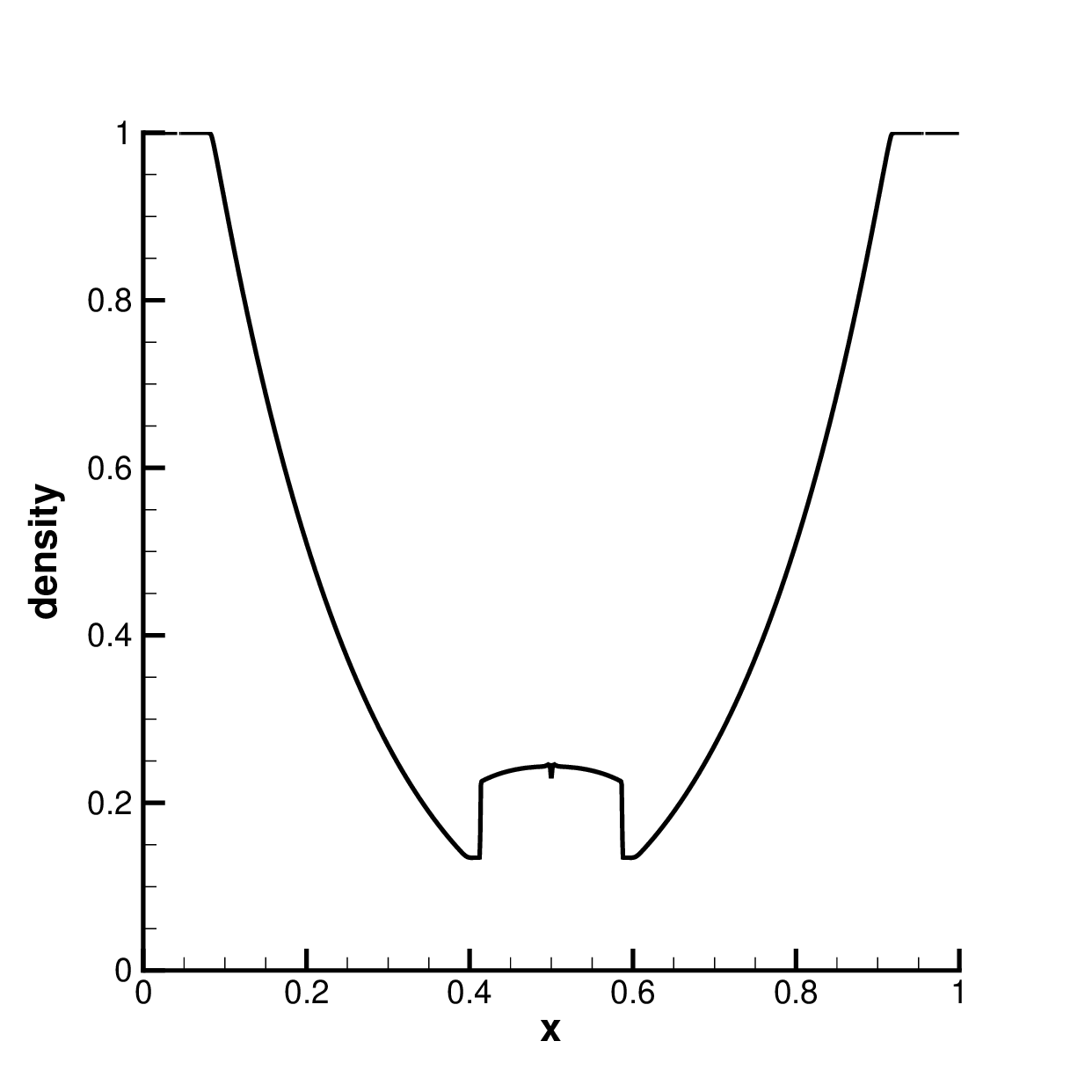}}
	\hspace{0mm}
	\caption{\footnotesize The density distribution for the first 2-D Riemann problem at $t=0.2$ (left, 30 uniform contours are drawn) and the section of the density along $y=x$ (right).}
	\label{fig:RP1}
\end{figure}

The initial conditions for the first Riemann problem are
\begin{equation}
	(\rho_0, u_0, v_0, p_0) = \begin{cases}
		\begin{aligned}
			&(1, 0.6233, 0.6233, 1.5), \quad\quad  &x>0.5,\, y>0.5,\\
			&(0.389, -0.6233, 0.6233, 0.4),\quad\quad  &x<0.5,\, y>0.5,\\
			&(1, -0.6233, -0.6233, 1.5), \quad\quad  &x<0.5,\, y<0.5,\\
			&(0.389, 0.6233, -0.6233, 0.4), \quad\quad  &x>0.5,\, y<0.5.
		\end{aligned}
	\end{cases}\notag
\end{equation}
The density counters and the density distribution along the line $y=x$ at $t=0.2$ are shown in \Cref{fig:RP1}. The low density region is well resolved, and the continuous transition from continuous flows to the presence of transonic shocks is well captured by the current scheme.

\begin{figure}[!htb]
	\centering
	\subfigure{	\includegraphics[width=0.45\linewidth]{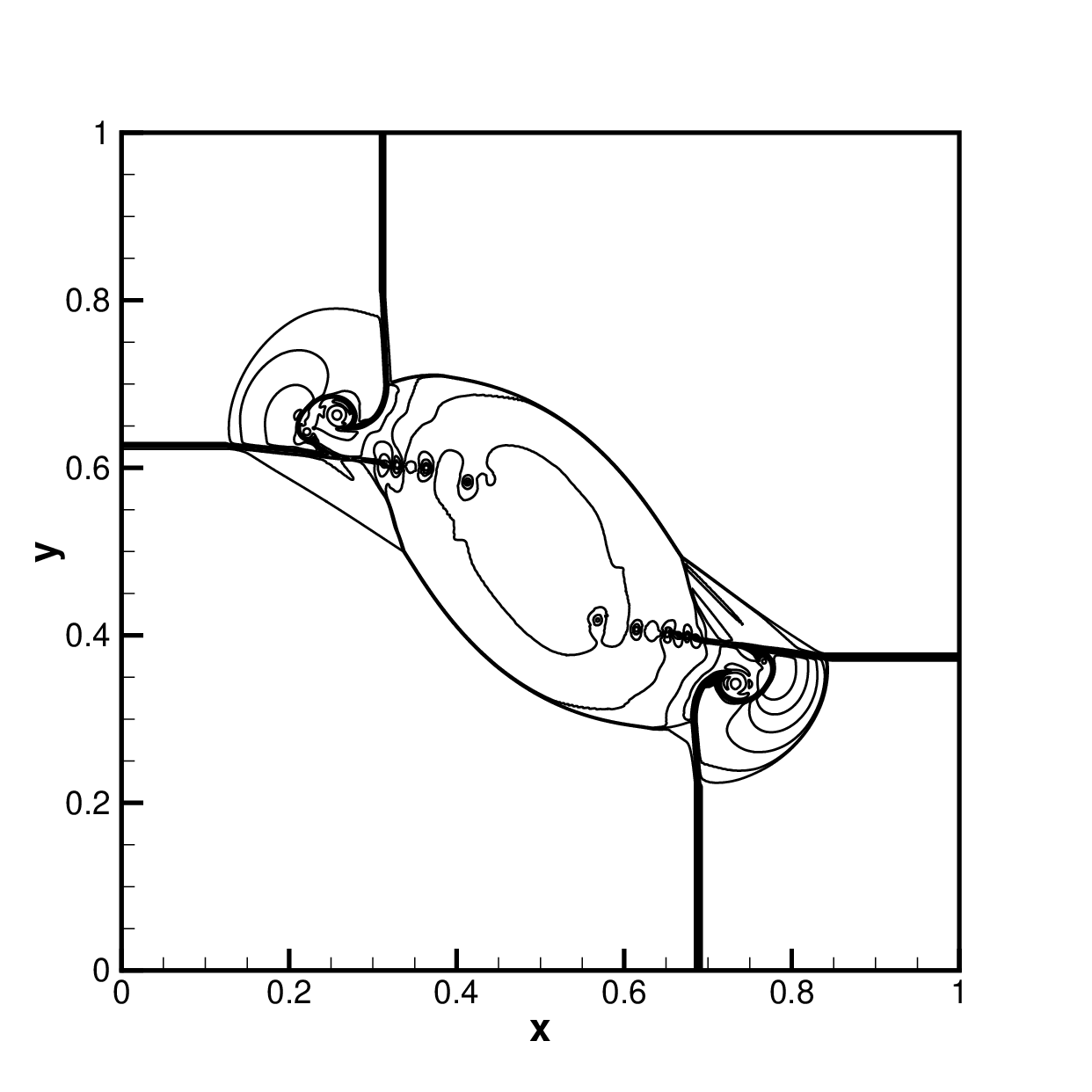}}
	\hspace{0mm}
	\subfigure{	\includegraphics[width=0.45\linewidth]{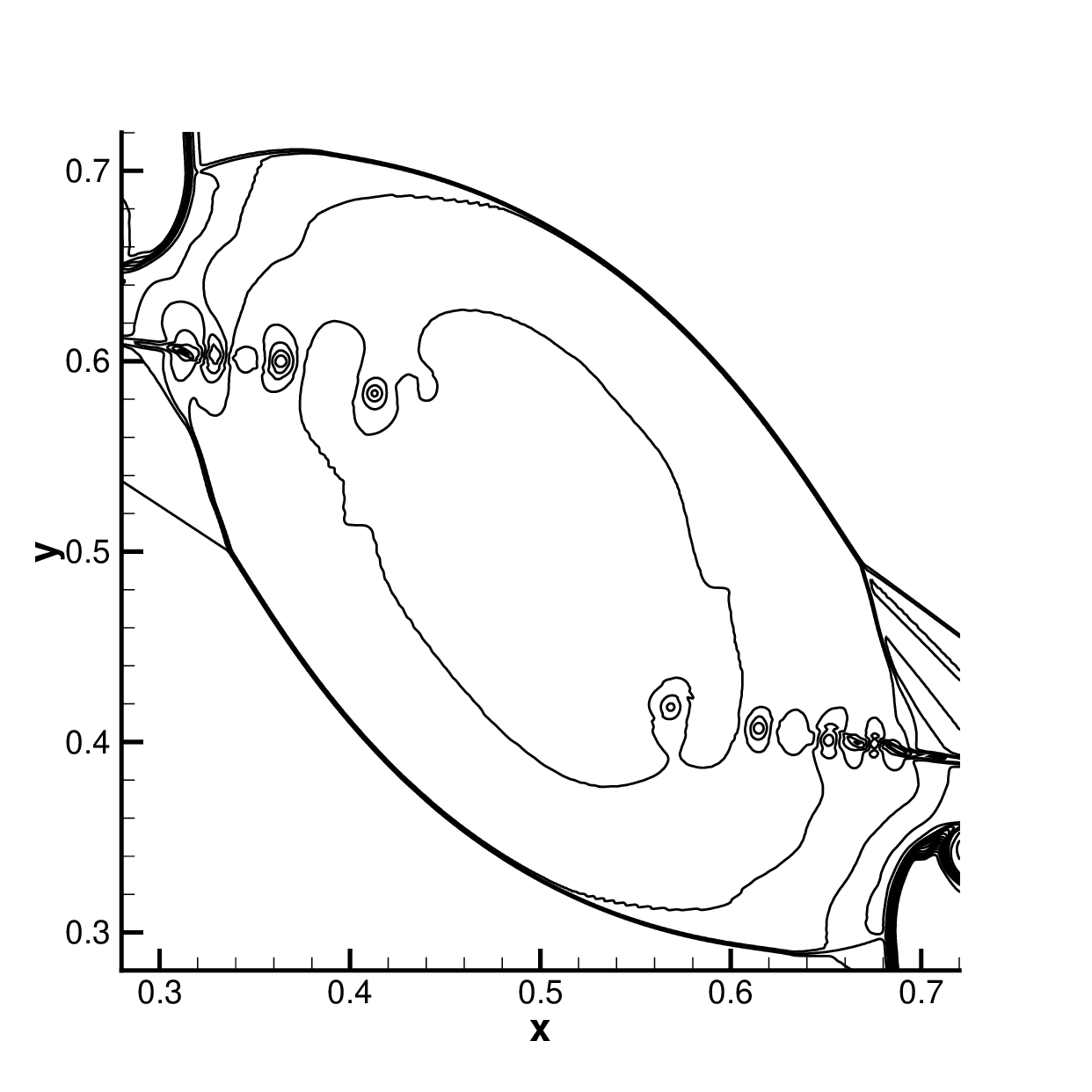}}
	\hspace{0mm}
	\caption{\footnotesize The density distribution for the second 2-D Riemann problem at $t=0.25$ (left) and the local enlargement (right). 20 uniform contours are drawn.}
	\label{fig:RP2}
\end{figure}

The initial conditions for the second Riemann problem are
\begin{equation}
	(\rho_0, u_0, v_0, p_0) = \begin{cases}
		\begin{aligned}
			&(1, -0.75, -0.5, 0.75), \quad\quad  &x>0.5,\, y>0.5,\\
			&(2, -0.75, 0.5, 0.75), \quad\quad  &x<0.5,\, y>0.5,\\		
			&(1, 0.75, 0.5, 0.75),\quad\quad  &x<0.5,\, y<0.5,\\
			&(3, 0.75, -0.5, 0.75), \quad\quad  &x>0.5,\, y<0.5.		
		\end{aligned}
	\end{cases}\notag
\end{equation}
The density counters at $t=0.25$ and the local enlargement are shown in \Cref{fig:RP2}. The small scaled vortices are well captured by the current scheme.


\subsubsection{Double Mach reflection problem}
The classical double Mach reflection problem is tested here to demonstrate the robustness of the current scheme in capturing strong shocks. This problem was initially proposed by Woodward and Colella \cite{WC1984}.

The computational domain is $[0, 4] \times [0, 1]$, and a solid wall lies at the bottom of the computational domain starting from $x = 1/6$. 
A right-moving shock initially makes a $\pi/3$ angle with the positive direction of the $x$-axis and hits the bottom at $(x, y) = (1/6, 0)$. The pre-shock state is $(\rho, u, v, p) = (1.4, 0, 0, 1)$ and the Mach number of the shock is $\Ma=10$.

The reflecting boundary condition is applied along the wall, while for the rest of the bottom boundary, the exact post-shock fluid state is imposed. At the top boundary, the flow variables are set to follow the motion of the Mach 10 shock. Uniform meshes with $960 \times 240$ and $1920\times 480$ cells are used.

\begin{figure}[!htb]
	\centering
	\subfigure{	\includegraphics[width=0.8\linewidth]{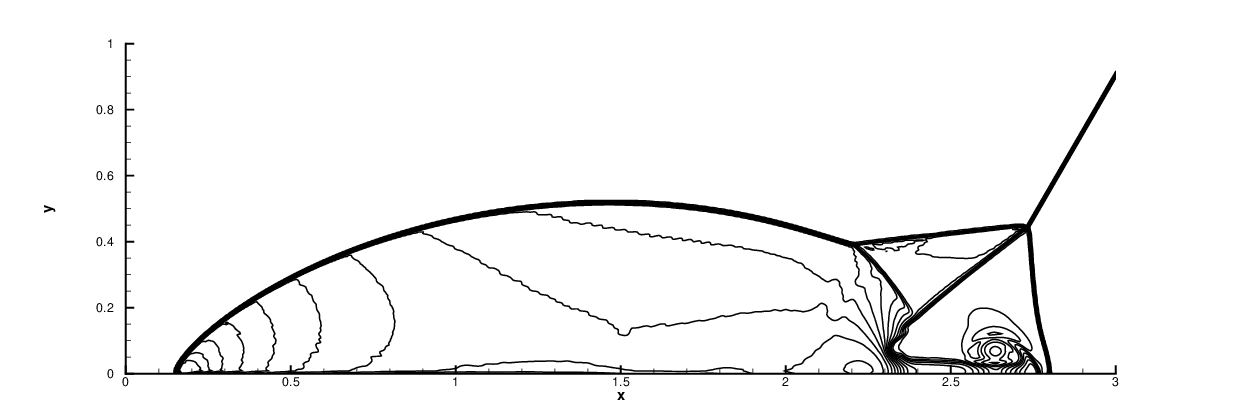}}
	\subfigure{	\includegraphics[width=0.8\linewidth]{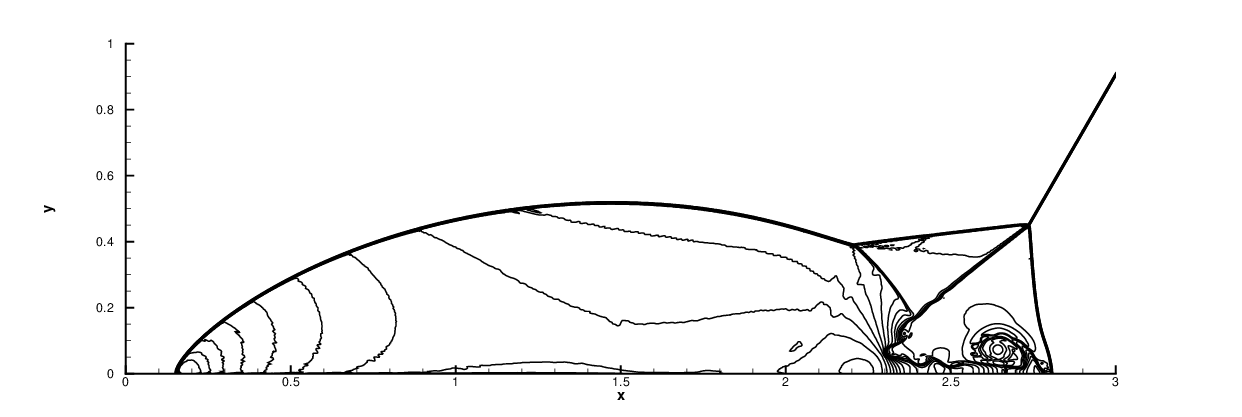}}
	\caption{\footnotesize Double Mach reflection problem: the density contours at $t=0.2$ with $960 \times 240$ (top) and $1920\times 480$ (bottom) cells. 30 uniform contours are displayed.}
	\label{fig:DoubleMach}
\end{figure}

\begin{figure}[!htb]
	\centering
	\subfigure{\includegraphics[width=0.45\linewidth]{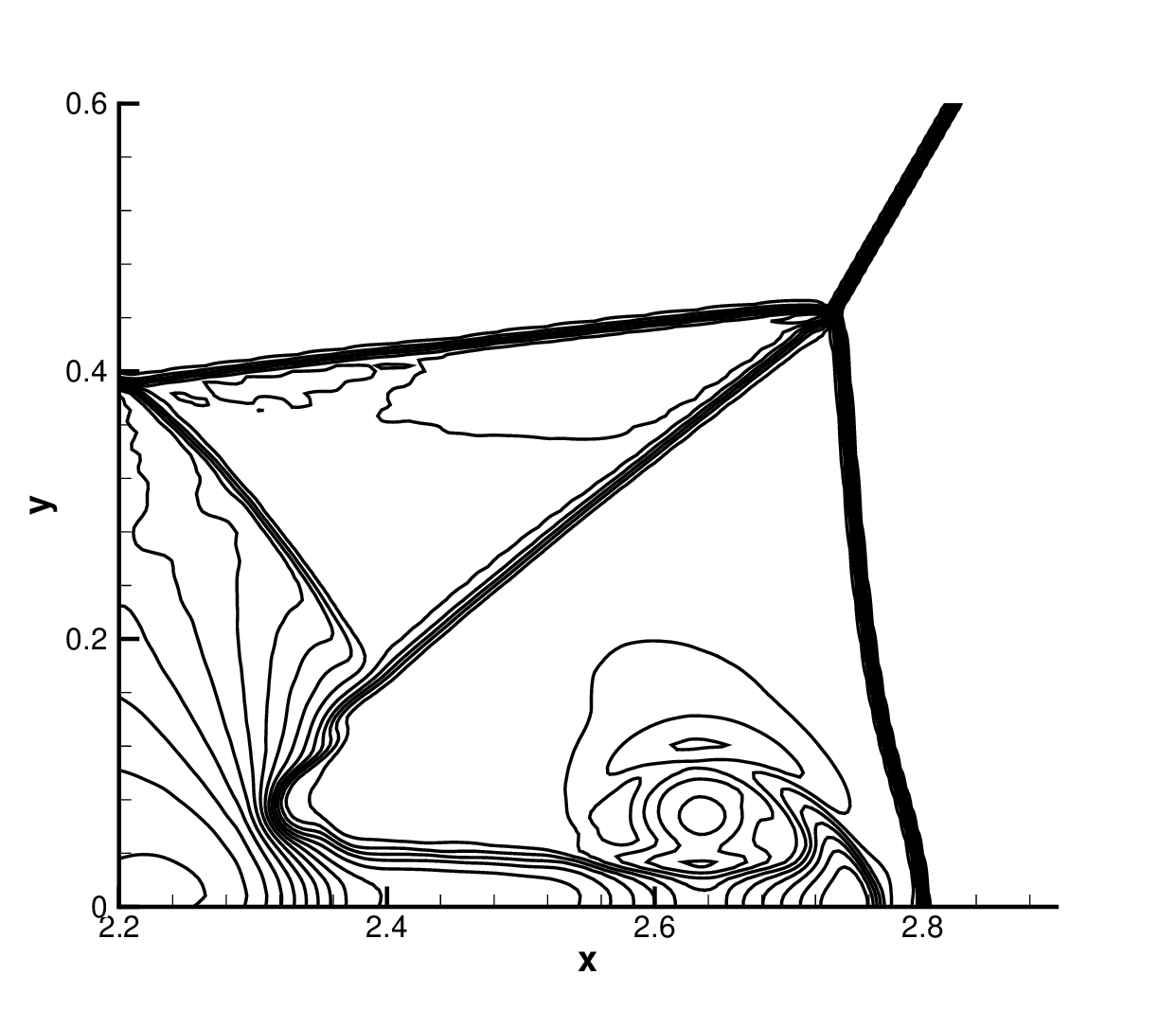}}	
	\subfigure{	\includegraphics[width=0.45\linewidth]{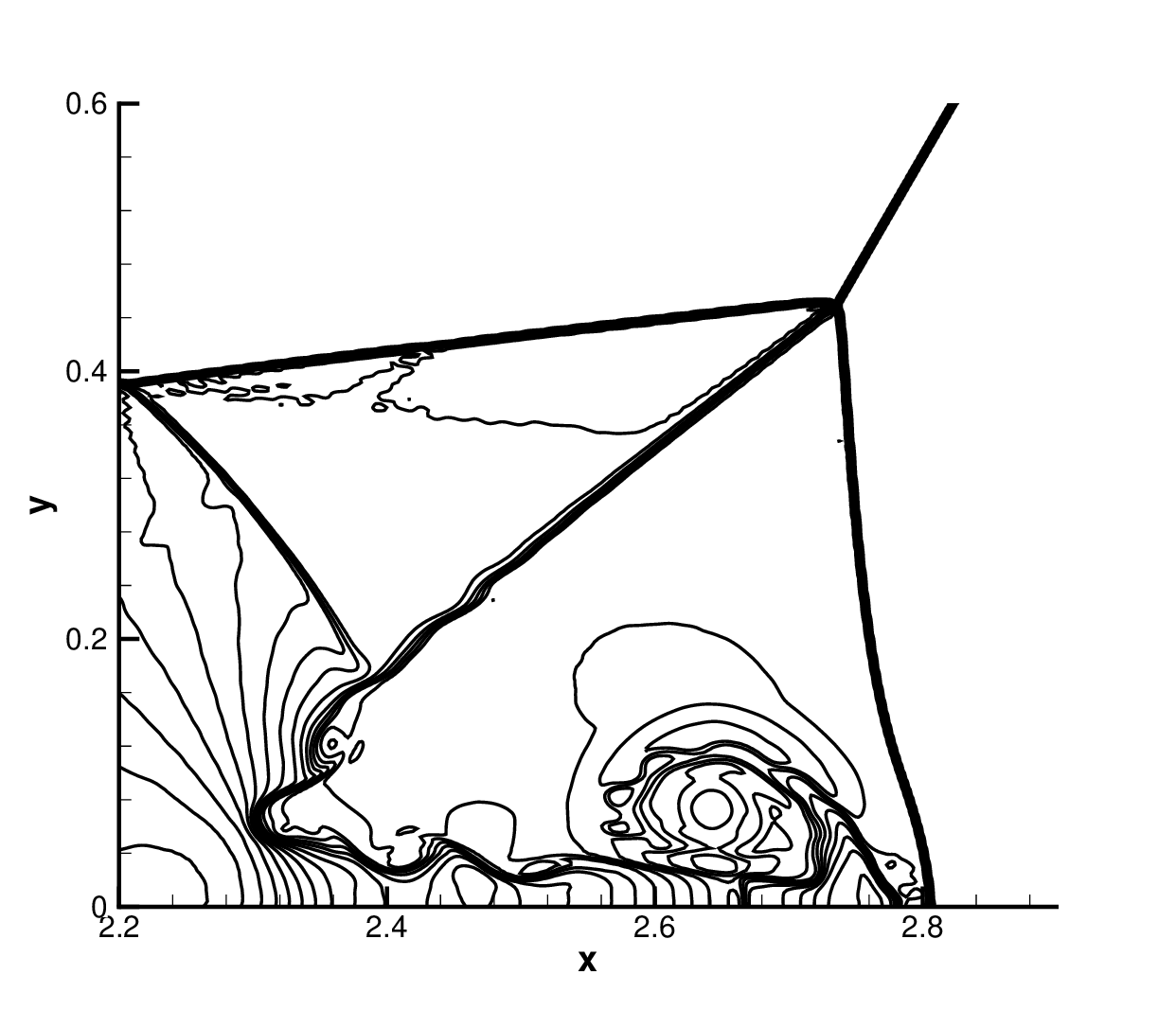}}
	\caption{\footnotesize Double Mach reflection problem: enlarged density contours around the triple point at $t=0.2$ with $960 \times 240$ (left) and $1920\times 480$ (right) cells. 30 uniform contours are displayed}.
	\label{fig:DoubleMach2}
\end{figure}

At $t=0.2$, the density contours and the local enlargements for the domain $[0,3]\times[0,1]$ are shown in \Cref{fig:DoubleMach} and \Cref{fig:DoubleMach2}. It is observed that the current scheme resolves the flow structure under the triple Mach stem clearly.

\subsection{Viscous flows}

\vspace{2mm}
\subsubsection{Couette flow problem}
The first test case of viscous flows is the Couette flow problem to verify the accuracy of the present numerical scheme for viscous smooth flows. The Couette flow is the laminar viscous flow between two parallel plates. The bottom plate with a fixed temperature $T_b$ is stationary, while the top one with a fixed temperature $T_1$ is moving at a constant horizontal speed $U$. The two plates are separated by a distance $H$. 

There exists an exact solution under the simplification that the viscosity coefficient $\mu$ is a constant and the Mach number is low enough to ensure nearly incompressible assumptions. The analytic steady-state solution is 
\begin{equation}\notag
	\begin{aligned}
		&u= \frac{y}{H}U,\quad v=0,\\
		&p=p_\infty,\quad \rho = \frac{p}{T},\\
		&T=T_0+\frac{y}{H}(T_1-T_b)+\frac{y}{H}(1.0-\frac{y}{H})\frac{\Pr U^2}{2C_p},
	\end{aligned}
\end{equation}
where the $p_\infty$ is the pressure at the upper wall, $\Pr$ is the Prandtl number, and $C_p$ is the specific heat capacity at constant pressure. 

In this test, we take $H=1$, $T_1=1$, $T_b=0.85$ and $p_\infty=1$. The Mach number for the upper wall is set as $\Ma=U/\sqrt{\gamma T_1}=0.1$. The density at the upper wall is $\rho_1=p_1/ T_1=1$ and the Reynolds number is defined as $\Re=\rho_1 U H/\mu=100$. The computational domain is a rectangle $[0,2H]\times[0,H]$. A uniform mesh with $2N\times N$ cells is used. All the boundary conditions are handled by ghost cells. The flow variables at ghost cells are assigned with the exact solution of the steady state. The initial value is taken as the small perturbation of the exact solution. The flow field is assumed to reach a steady-state when the $L^2$-norm of the temperature residual is less than $10^{-14}$.

The errors and convergence orders of the velocity and the temperature are presented in \Cref{Tab:CouetteF-vel} and \Cref{Tab:CouetteF-temp}. It is observed that the designed order is achieved for viscous flows.

\begin{table}[!htbp]
	\centering
	\caption{\footnotesize Accuracy test for the Couette flow problem: the $L^1$- and $L^\infty$-errors and convergence orders for the velocity $u$.}   \label{Tab:CouetteF-vel}   
	\setlength{\tabcolsep}{3.5mm}
	\begin{tabular}{ccccc}
		\toprule
		
		N  &$L^1$ error& order  
		&$L^\infty$ error  & order  
		\\
		\hline
		10  &5.842e-07 &    &1.664e-06  &    
		\\
		20  &1.613e-07 &1.856  &4.487e-07   &1.891   
		\\
		40  &4.021e-08 &2.005   &1.218e-07   &1.881     
		\\
		80   &1.000e-08 &2.007   &3.138e-08    &1.957    
		\\
		160  &2.484e-09 &2.009   &7.894e-09   &1.991   
		\\
		\bottomrule
	\end{tabular}   
\end{table}  
\begin{table}[!htbp]
	\centering
	\caption{\footnotesize Accuracy test for the Couette flow problem: the $L^1$- and $L^\infty$-errors and convergence orders for the temperature $T$.}   \label{Tab:CouetteF-temp}   
	\setlength{\tabcolsep}{3.5mm}
	\begin{tabular}{ccccc}
		\toprule
		
		N  &$L^1$ error& order  
		&$L^\infty$ error  & order  
		\\
		\hline
		10  &3.354e-06 &    &1.112e-05   &    
		\\
		20  &7.203e-07 &2.219  &2.014e-06   &2.464   
		\\
		40  &1.657e-07 &2.120   &4.766e-07   &2.079    
		\\
		80  &3.907e-08 &2.084   &1.143e-07   &2.059    
		\\
		160  &9.127e-09 &2.097   &2.701e-08   &2.082  
		\\
		\bottomrule
	\end{tabular}   
\end{table}

\vspace{2mm}
\subsubsection{Laminar boundary layer}
\begin{figure}[!htb]
	\centering
	\setlength{\abovecaptionskip}{0.0cm}
	\setlength{\belowcaptionskip}{0.0cm}
	\subfigure{	\includegraphics[width=0.7\linewidth]{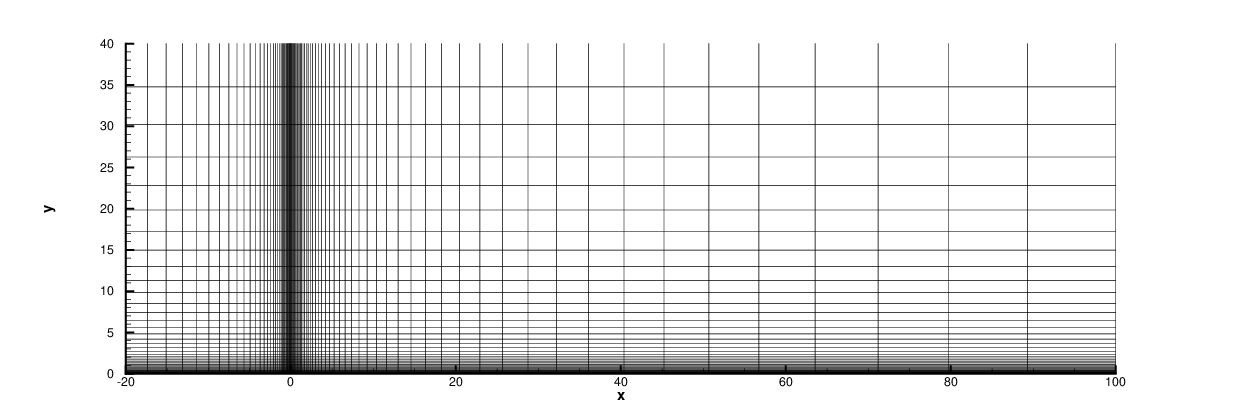}}
	\caption{\footnotesize Laminar boundary layer: the rectangular mesh.}
	\label{fig:mesh}
\end{figure}

A laminar boundary layer over a flat plate with the length $L=100$ is considered here. The Mach number of the inflow free-stream is $\Ma=0.15$ and the Reynolds number is $\Re=\frac{u_{\infty}L}{\mu}=10^5$, where $u_{\infty}$ is the horizontal velocity of the free-stream and $\mu$ is the dynamic viscosity coefficient. Note that for this free-stream Mach number, the flow is nearly incompressible so the incompressible laminar boundary layer theory can be used to assess the performance of our method. Similar tests can be found in \cite{DY2022,PXLL2016}. As shown in \Cref{fig:mesh}, the computational domain is set as $\Omega=[-20,100]\times[0,40]$, and the rectangular mesh consists of $86\times43$ cells. At the bottom boundary, the flat plate starts at $(x,y)=(0,0)$ and ends at $(x,y)=(100,0)$. There are $52$ cells allocated on the flat plate in the $x$-direction. The minimum mesh size is $0.04442$. At $y=0$, the non-slip and adiabatic boundary condition is applied for $0<x<100$, while a symmetric boundary condition is applied for $-20<x<0$. The inflow boundary condition is adopted for the inlet boundary $x=-20$. Moreover, the outflow boundary condition is enforced at the top boundary $y=40$ and exit boundary $x=100$. The initial state is given by the free-stream. When the steady state is achieved, the $L^2$-norm of the velocity residual is less than about $8\times10^{-8}$.

\begin{figure}[!htb]
	\centering
	\setlength{\abovecaptionskip}{0.0cm}
	\setlength{\belowcaptionskip}{0.0cm}
	\subfigure{	\includegraphics[width=0.45\linewidth]{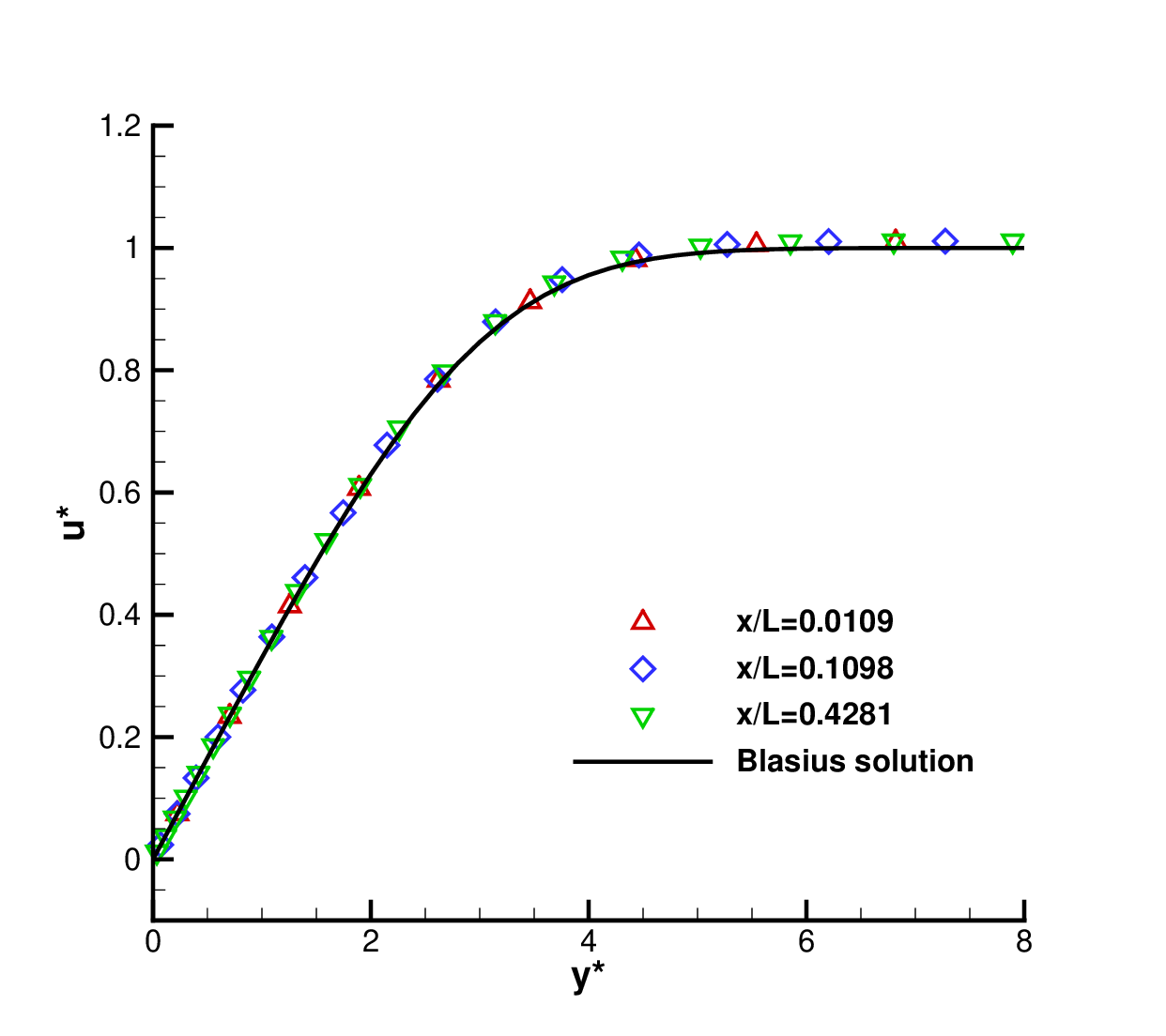}}
	\hspace{1mm}
	\subfigure{	\includegraphics[width=0.45\linewidth]{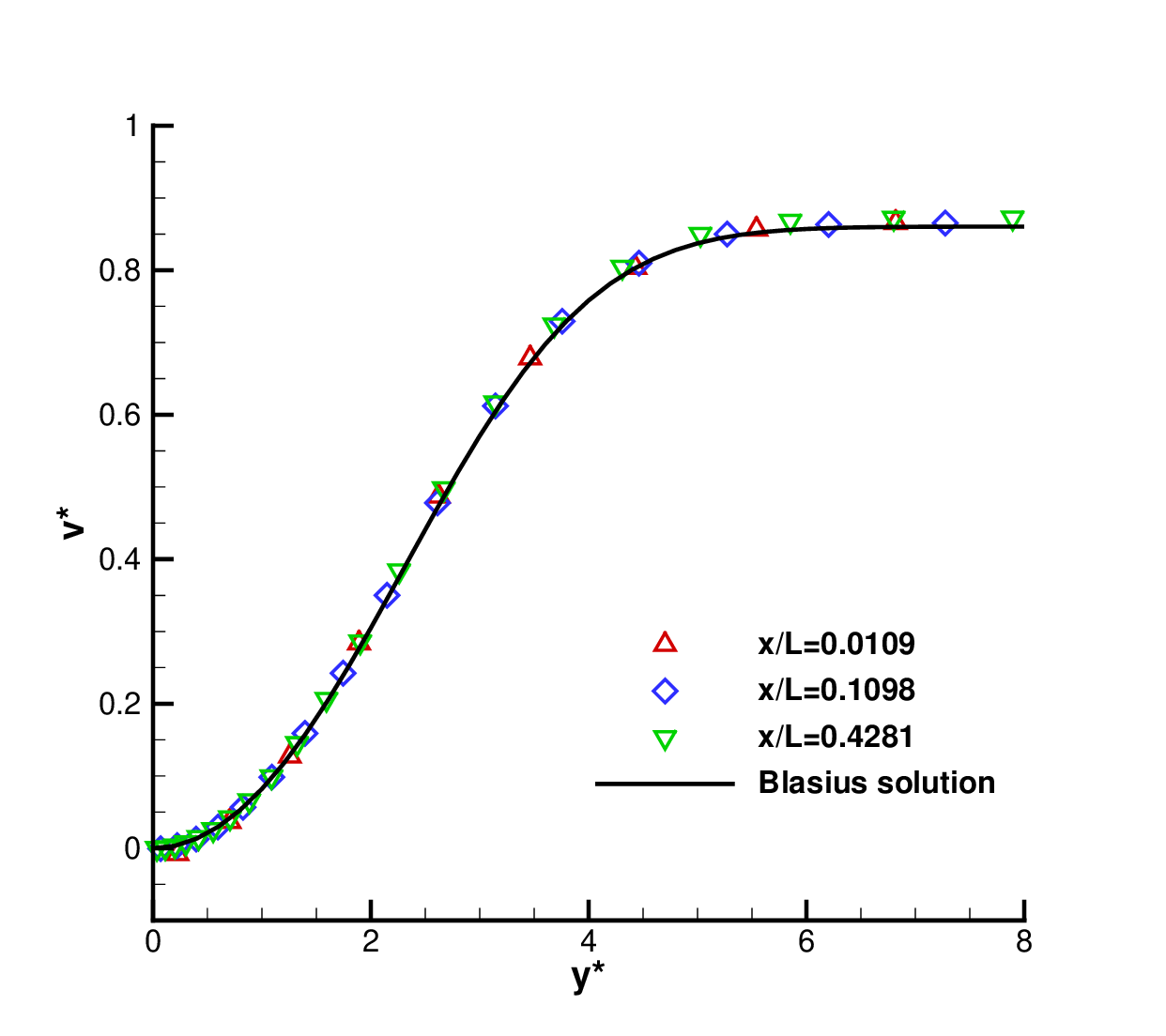}}
	\caption{\footnotesize Laminar boundary layer: streamwise velocity profiles at different locations for $\Re=10^5$.}
	\label{fig:Laminar}
\end{figure}

The streamwise velocity profiles at different locations are presented in \Cref{fig:Laminar}, where $u^*=\frac{u}{u_{\infty}}$, $v^*=\frac{v\sqrt{\Re_x}}{u_{\infty}}$, $y^*=\frac{y}{x}\sqrt{\Re_x}$, $\Re_x=\Re\frac{x}{L}$, and the solid lines are the exact Blasius solutions of $u^*$ and $v^*$. In all locations, the numerical solutions match the exact Blasius solution well. At the location close to the leading edge, the boundary layer profile can be accurately captured with only seven cells within the layer.


\vspace{2mm}

\subsubsection{Lid-driven cavity flow}
In order to further test the current scheme in capturing vortex flows, we test the lid-driven cavity flow, which is a typical benchmark for 2-D incompressible or low Mach number compressible viscous flows. The fluid is bounded by a unit square $\Omega=[0, 1]\times[0, 1]$. The upper wall moves with a horizontal speed $U=1$. Other walls are fixed. The non-slip and isothermal boundary conditions are applied to all boundaries with the temperature $T_{b}$. The Mach number is set as $\Ma=U/\sqrt{\gamma T_b} =0.15$. The initial flow is stationary with the density $\rho_0=1$ and the temperature $T_0=T_b$. Numerical simulations are conducted for three Reynolds numbers, i.e., $\Re = 400$, $1000$ and $3200$. A uniform mesh with $97\times97$ cells is used.

\begin{figure}[!htb]
	\centering
	\subfigure{	\includegraphics[width=0.45\linewidth]{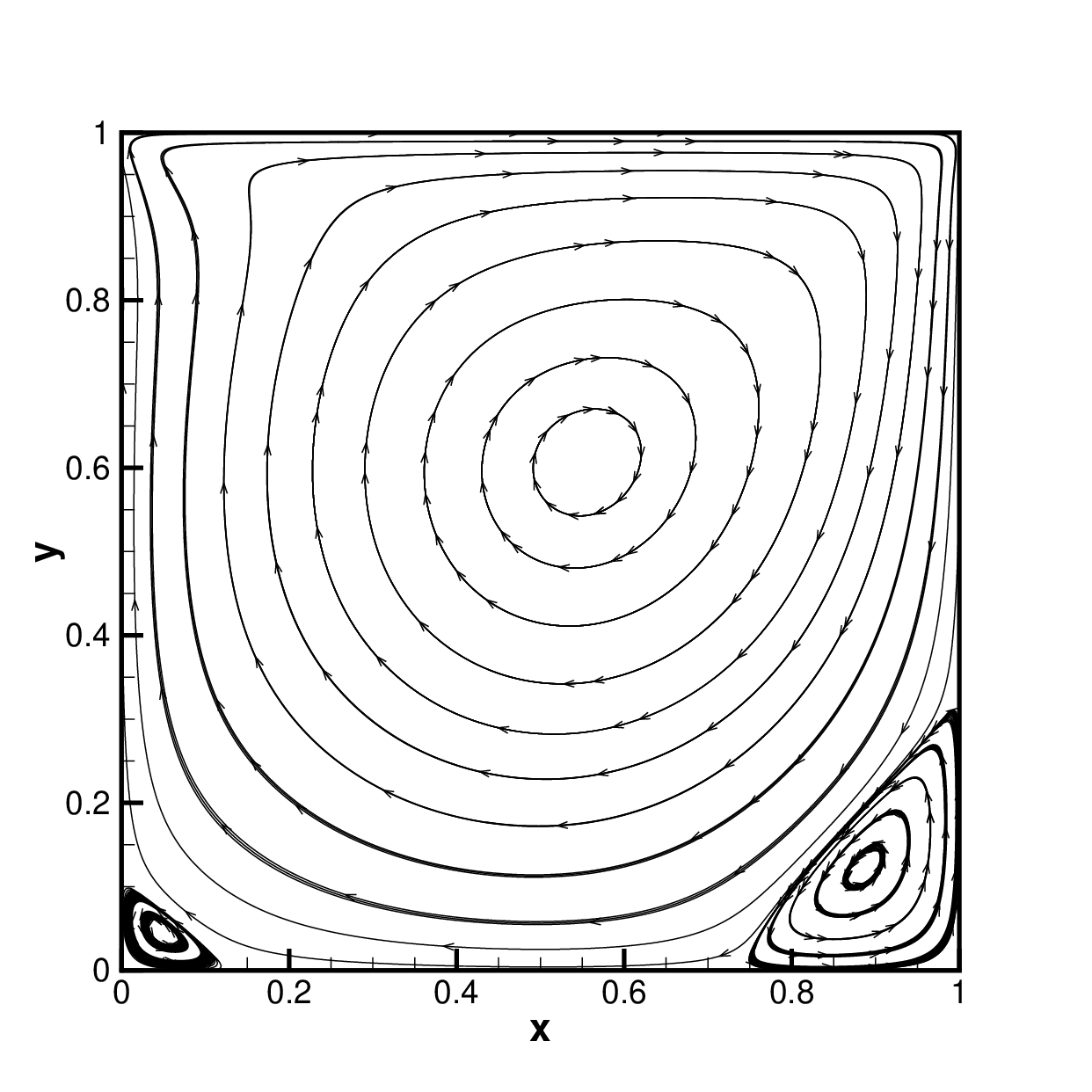}}
	\subfigure{	\includegraphics[width=0.45\linewidth]{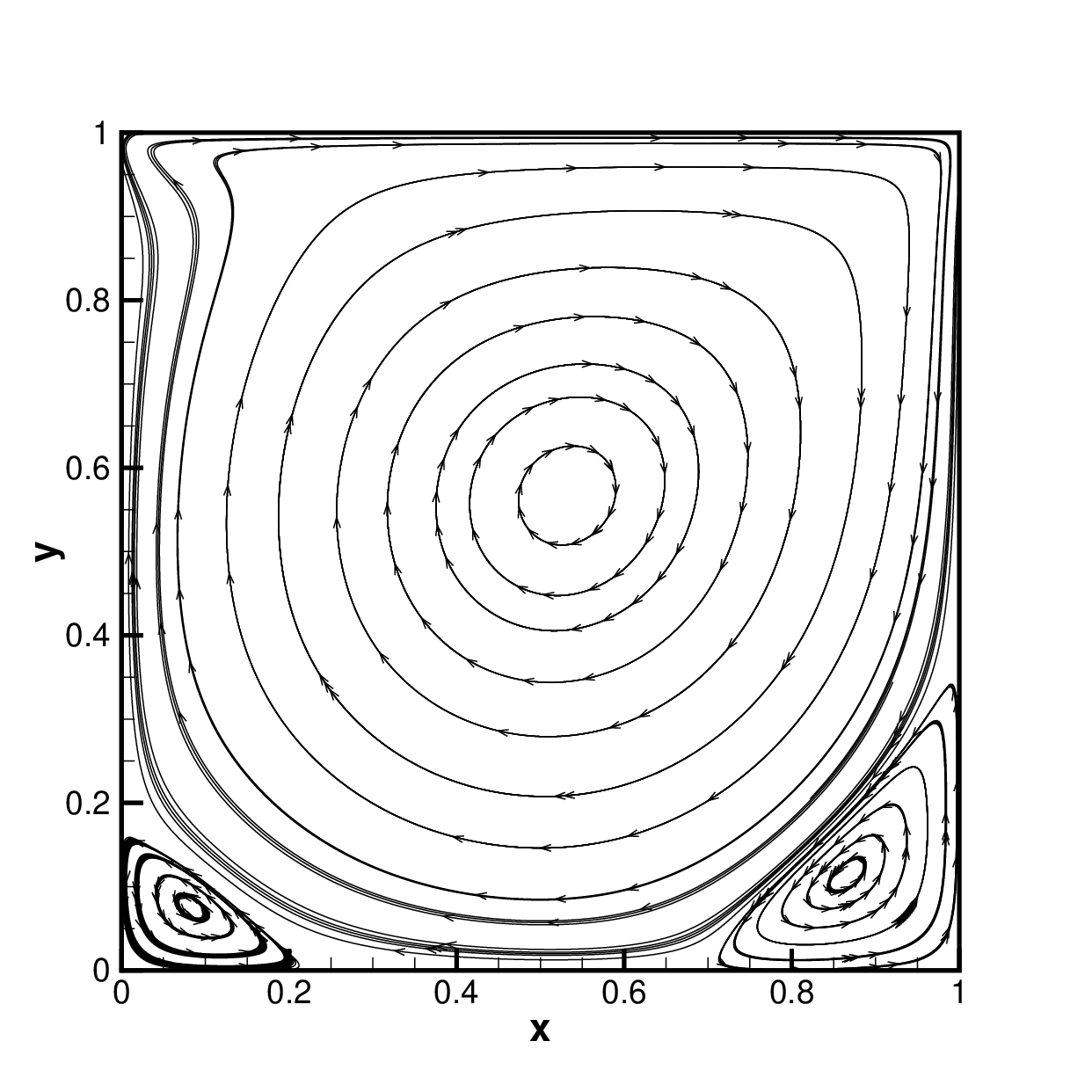}}
	\hspace{5mm}
	\subfigure{	\includegraphics[width=0.45\linewidth]{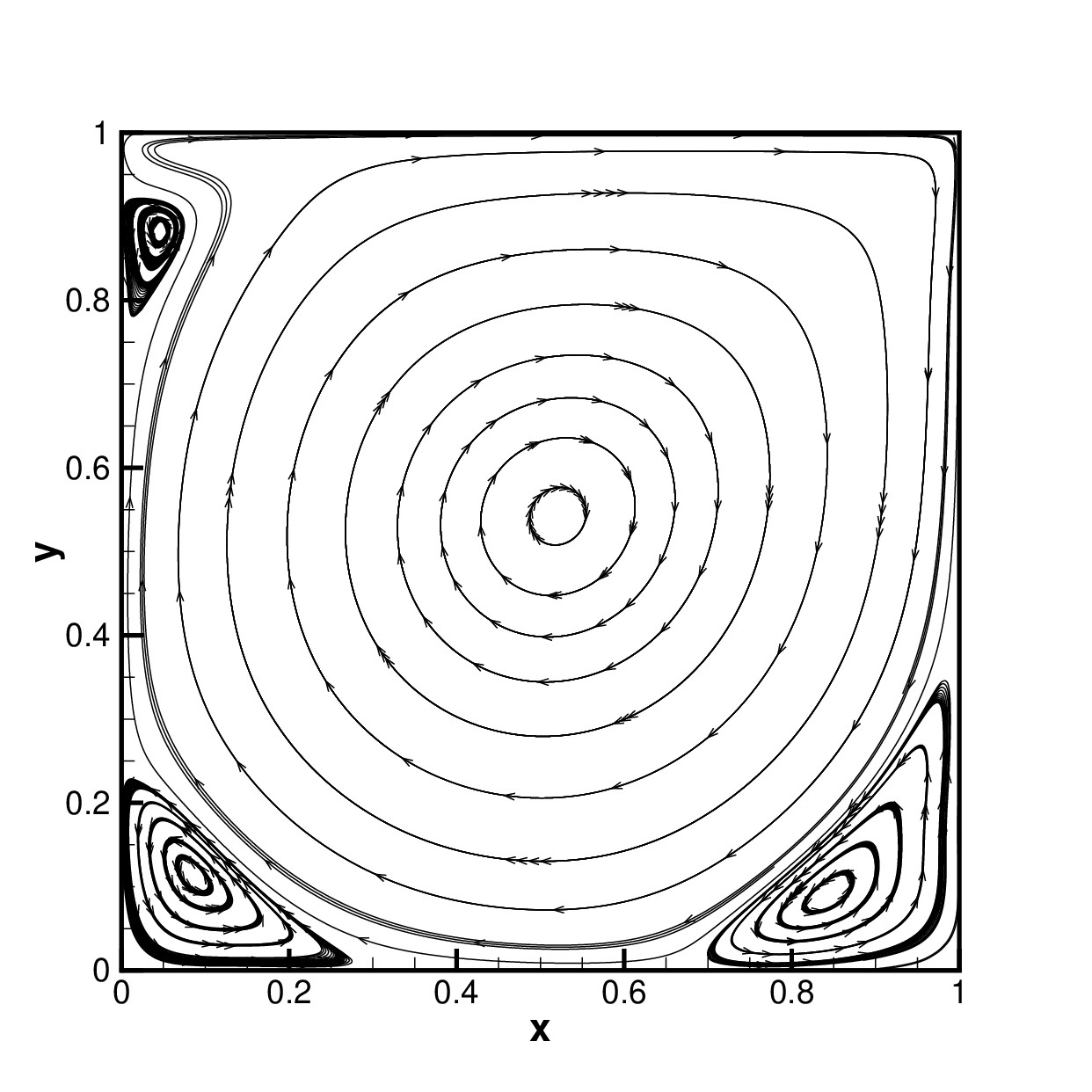}}
	\caption{\footnotesize Lid-driven cavity flow: the streamlines for $\Re=400$ (upper-left), $1000$ (upper-right) and $3200$ (lower) at $t=150$ with $97\times97$ cells.}
	\label{fig:meshandstreamlines}
\end{figure}

\begin{figure}[!htb]
	\centering
	\footnotesize
	\subfigure{	\includegraphics[width=0.45\linewidth]{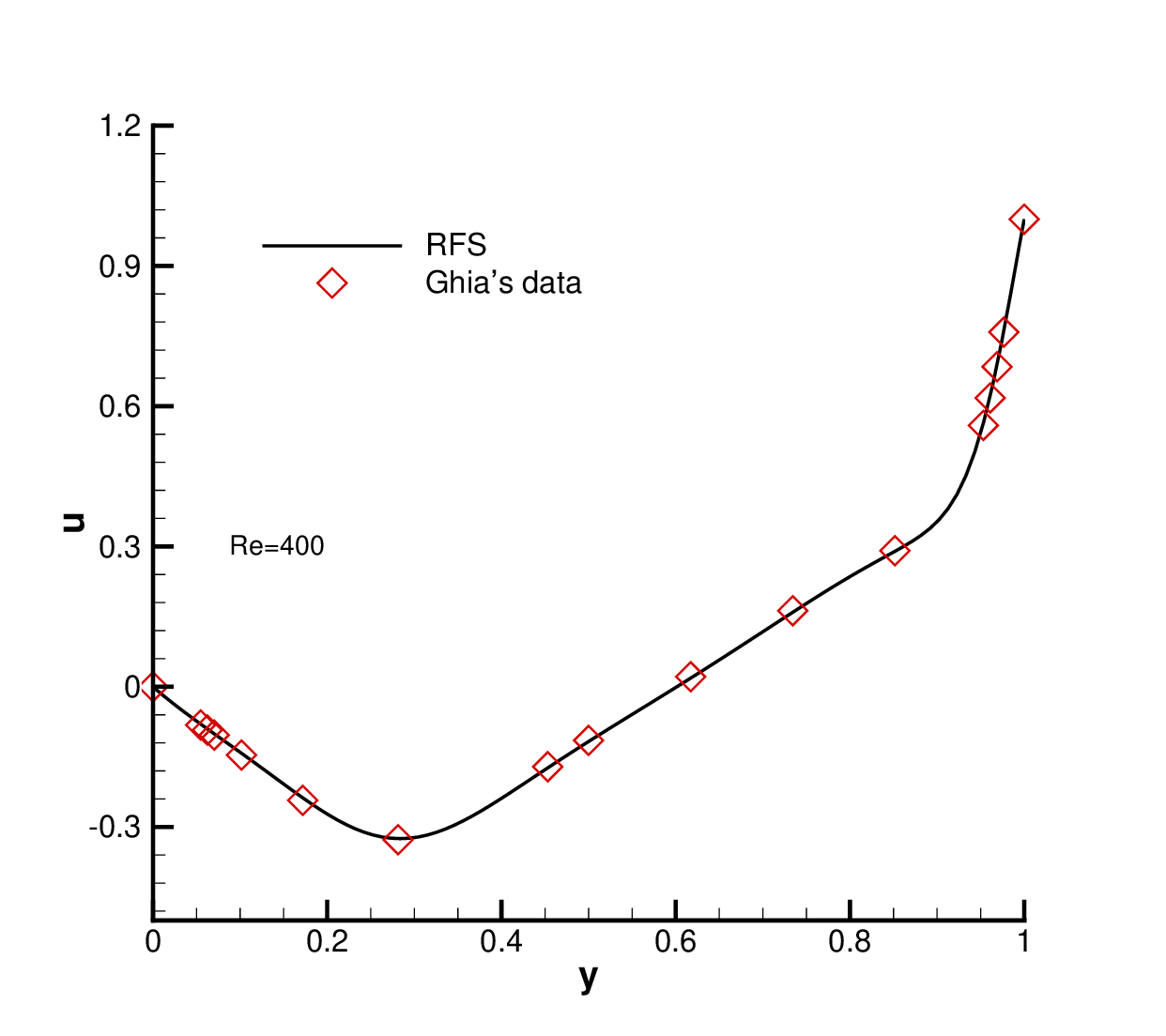}}
	\hspace{5mm}
	\subfigure{	\includegraphics[width=0.45\linewidth]{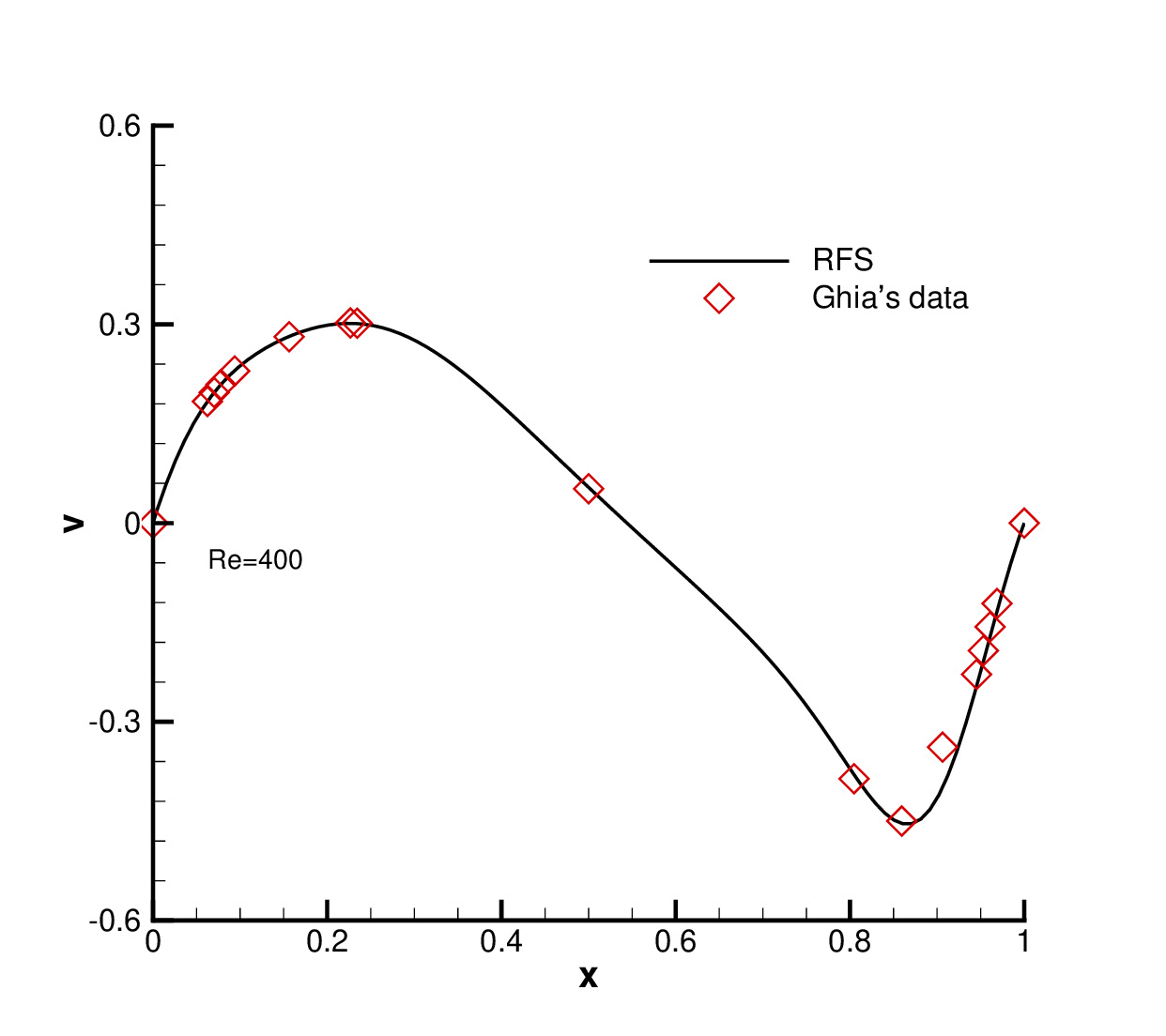}}
	\subfigure{	\includegraphics[width=0.45\linewidth]{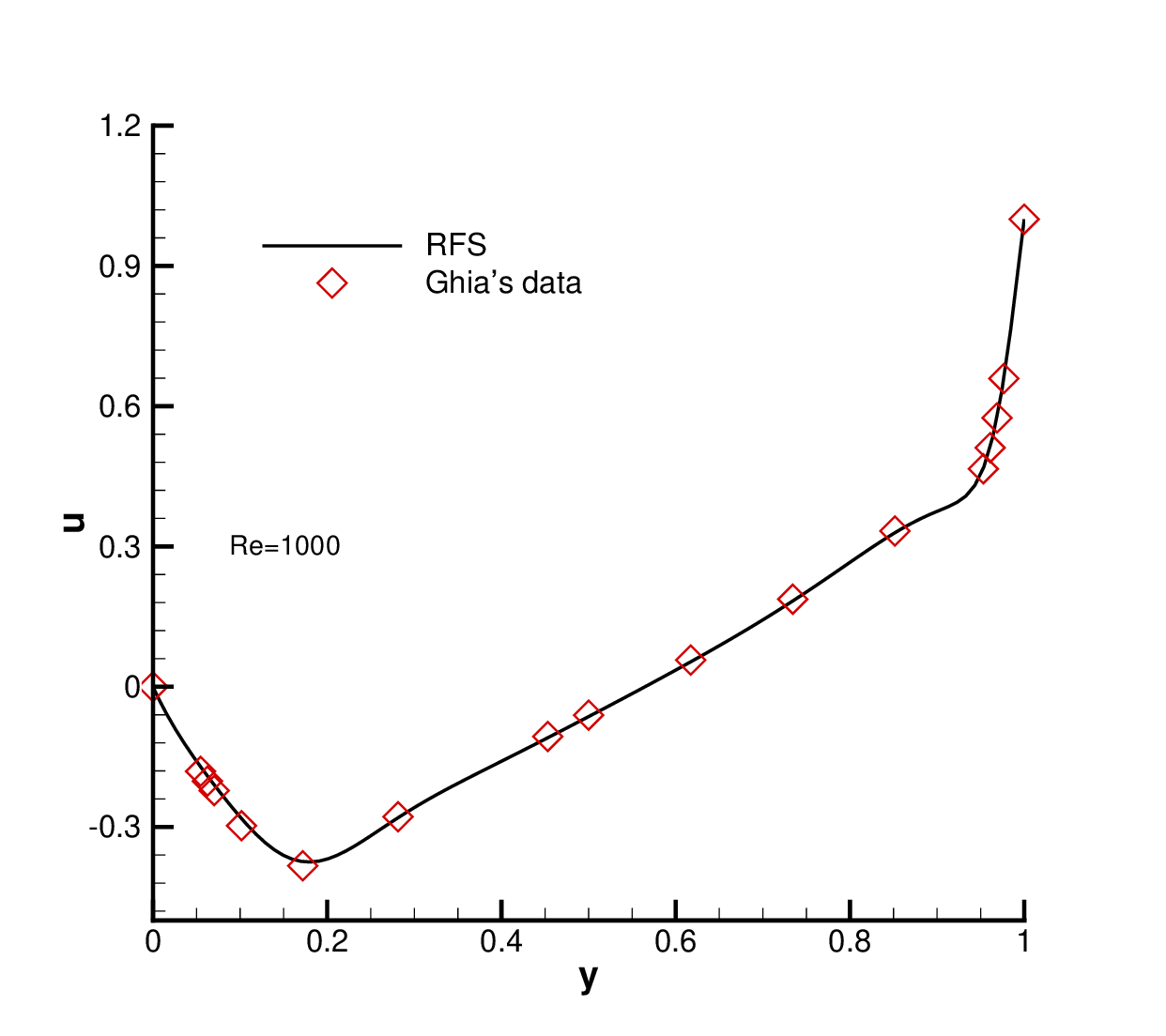}}
	\hspace{5mm}
	\subfigure{	\includegraphics[width=0.45\linewidth]{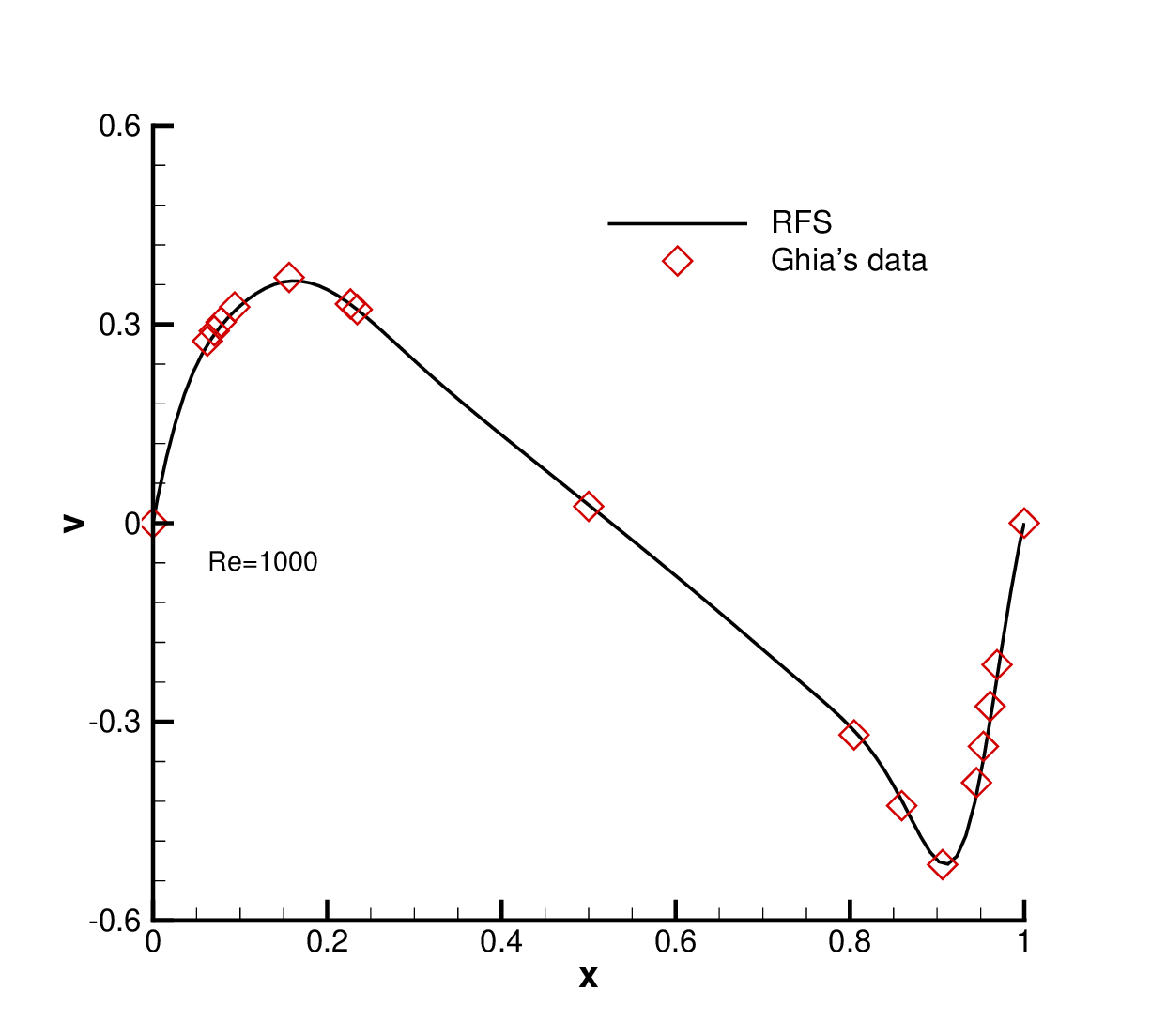}}
	\subfigure{	\includegraphics[width=0.45\linewidth]{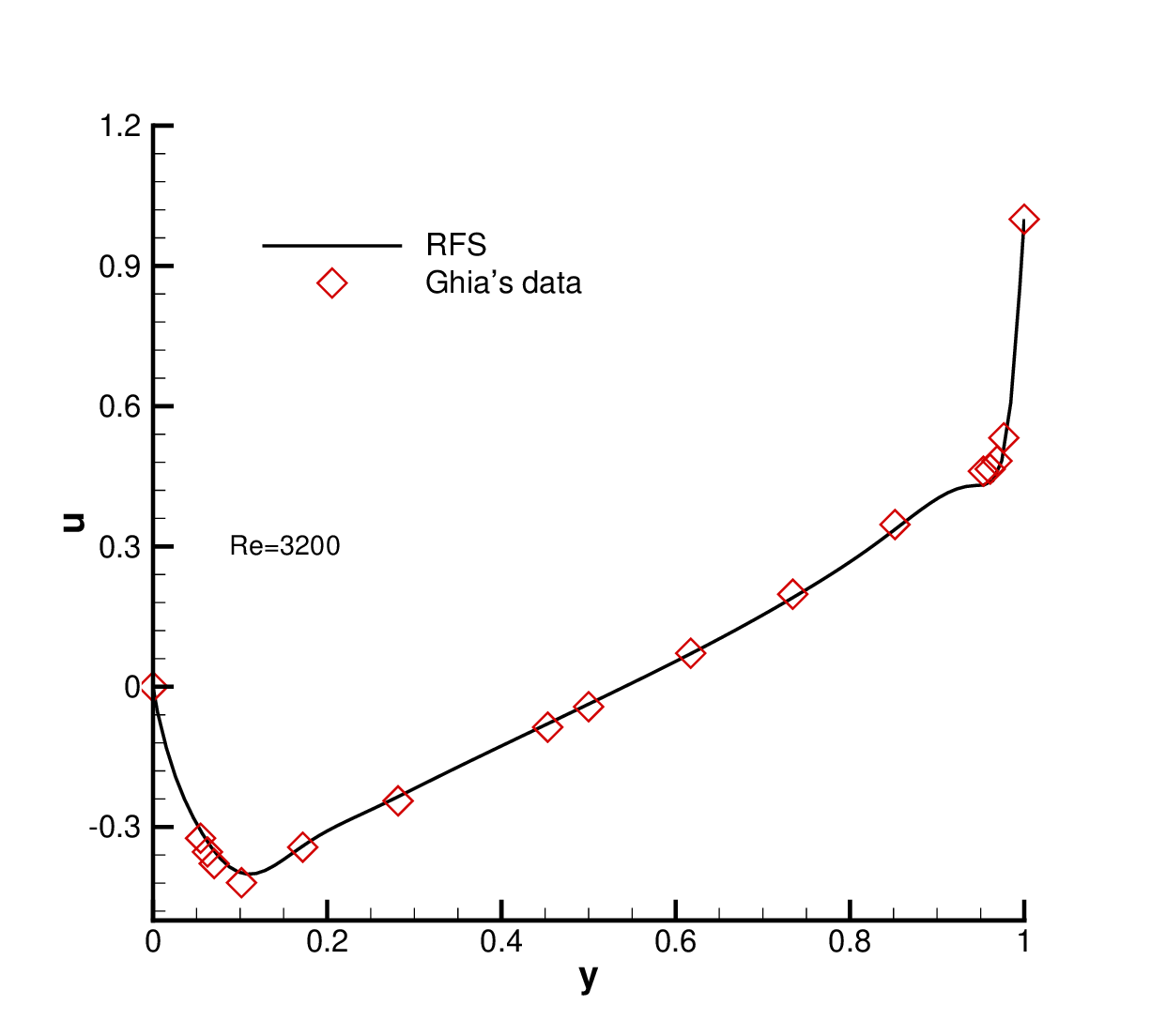}}
	\hspace{5mm}
	\subfigure{	\includegraphics[width=0.45\linewidth]{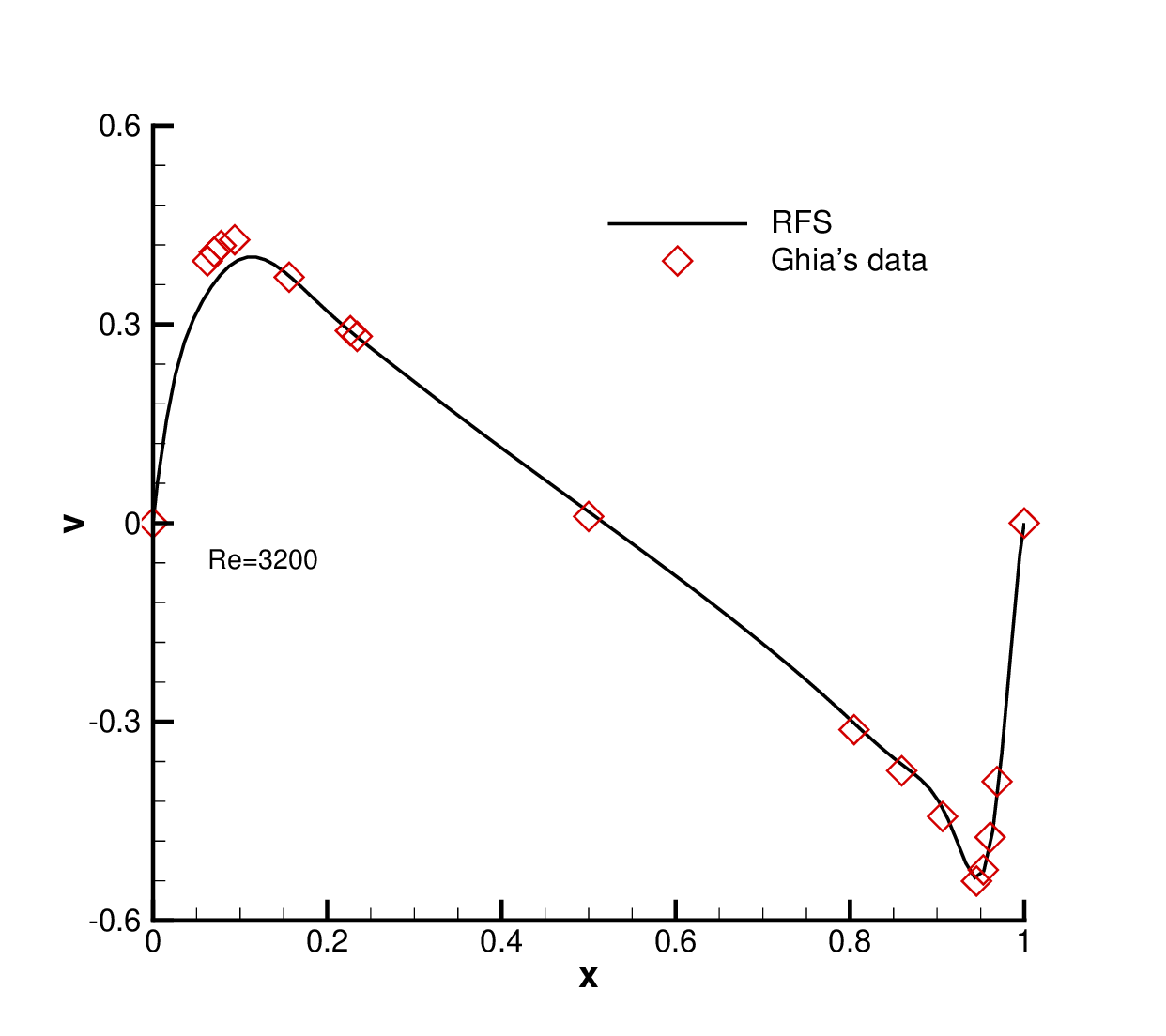}}
	\caption{\footnotesize Lid-driven cavity flow: $u$-velocity along vertical center line left column and $v$-velocity along horizontal center line right column for $\Re=400$ (top), $1000$ (middle) and $3200$ (bottom) at $t=150$ with $97\times97$ cells. The reference data are from Ghia et al.~\cite{GGS1982}.}
	\label{fig:uv}
\end{figure}

The streamlines at $t=150$ are shown in \Cref{fig:meshandstreamlines}. The primary flow structures including the primary and secondary vortices are well captured. The results of $u$-velocities along the center vertical line and $v$-velocities along the center horizontal line are shown in \Cref{fig:uv}. The simulation results match well with the benchmark data \cite{GGS1982}. The results demonstrate the high-resolution property of the current scheme.

\vspace{2mm}

\subsubsection{2-D viscous shock tube problem}
To demonstrate the performance of the current scheme in simulating high-speed viscous flows, we test the viscous shock tube problem, which has been studied extensively \cite{DT2009}. The flow is bounded by a unit square and is driven by a strong initial discontinuity at the center. Complex unsteady flow structures emerge because of interactions between the incident boundary layer and the reflected shock wave. In this case, an ideal gas is at rest in a unit square $[0,1]\times [0,1]$, and a membrane located at $x=0.5$ separates two different states of the gas. The dimensionless initial states are
\begin{equation}\notag
	(\rho, u, v, p) = \begin{cases}
		\begin{aligned}
			&(120, 0, 0, 120/\gamma ), \,\,\,\,\,\,\, &0\leq x\leq 0.5,\\
			&(1.2, 0, 0, 1.2/\gamma ), \,\,\,\,\,\,\, &0.5\leq x\leq 1.0.\\		
		\end{aligned}
	\end{cases}
\end{equation}
The membrane is removed at time zero and wave interaction occurs. A shock wave, followed by a contact discontinuity, moves right with a Mach number $\Ma=2.37$ and reflects at the right end wall. After the reflection, it interacts with the contact discontinuity. The contact discontinuity and shock wave interact with the horizontal wall and create a thin boundary layer during their propagation. The solution will develop complex two-dimensional shock/shear/boundary-layer interactions and requires not only strong robustness but also high resolution of a numerical scheme.

This case is tested in the computational domain $[0, 1] \times[0, 0.5]$. The symmetric boundary condition is applied on the upper boundary, and the non-slip and adiabatic boundary conditions are adopted on other boundaries. The Reynolds number $\Re=200$ and $1000$ are chosen, which are based on constant dynamic viscosity coefficients $\mu=0.005$ and $0.001$, respectively. 

\begin{figure}[!htb]
	\centering
	\subfigure{	\includegraphics[width=0.9\linewidth]{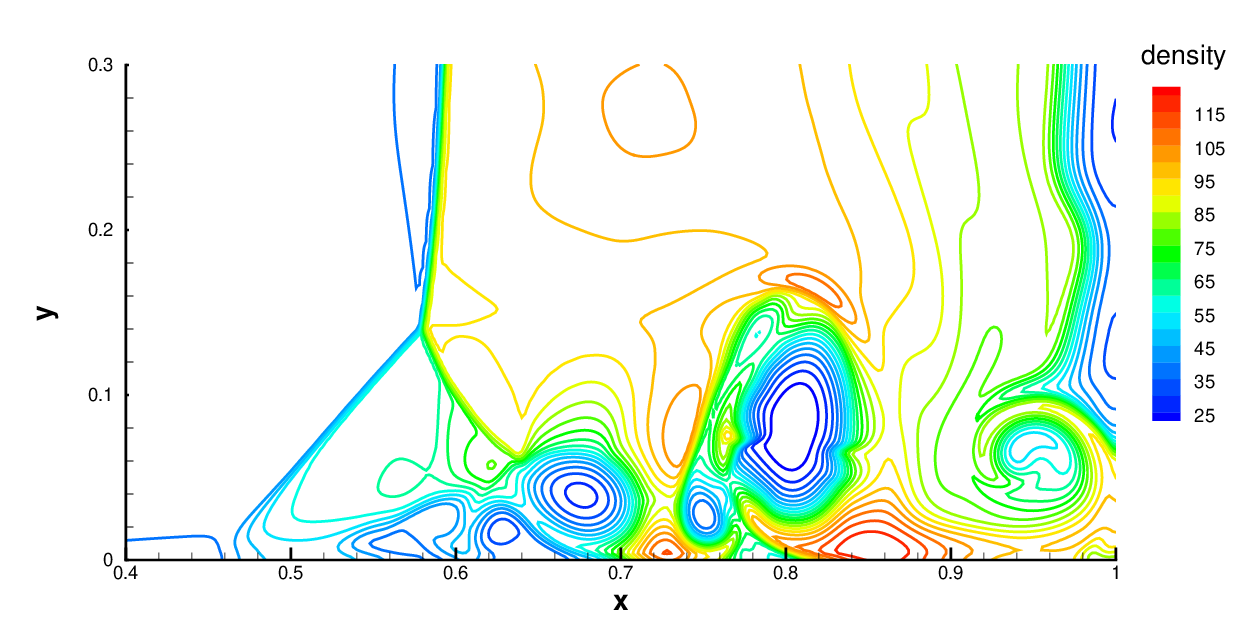}}
	\subfigure{	\includegraphics[width=0.9\linewidth]{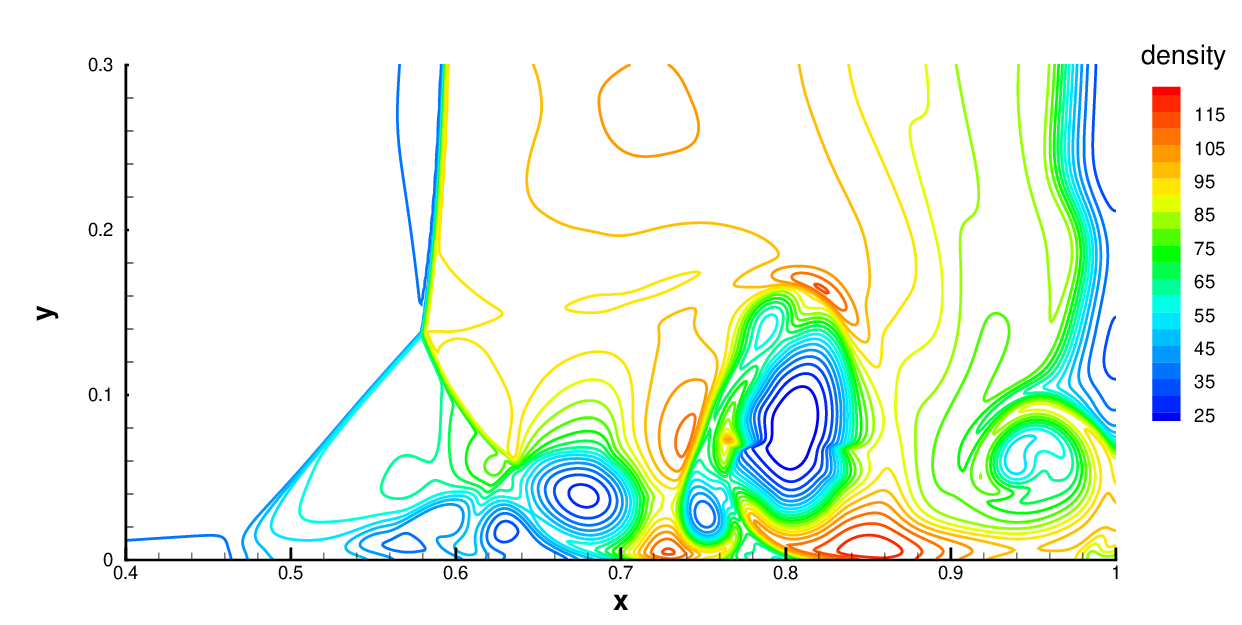}}
	\caption{\footnotesize Reflected shock-boundary layer interaction: density contours at $t=1.0$ for $\Re=200$ with the mesh size $h=1/500$ (top) and $1/1000$ (bottom). 20 uniform contours from 25 to 120 are drawn.}
	\label{fig:VSTPRe200}
\end{figure}

\begin{figure}[!htb]
	\centering
	\includegraphics[width=0.9\linewidth]{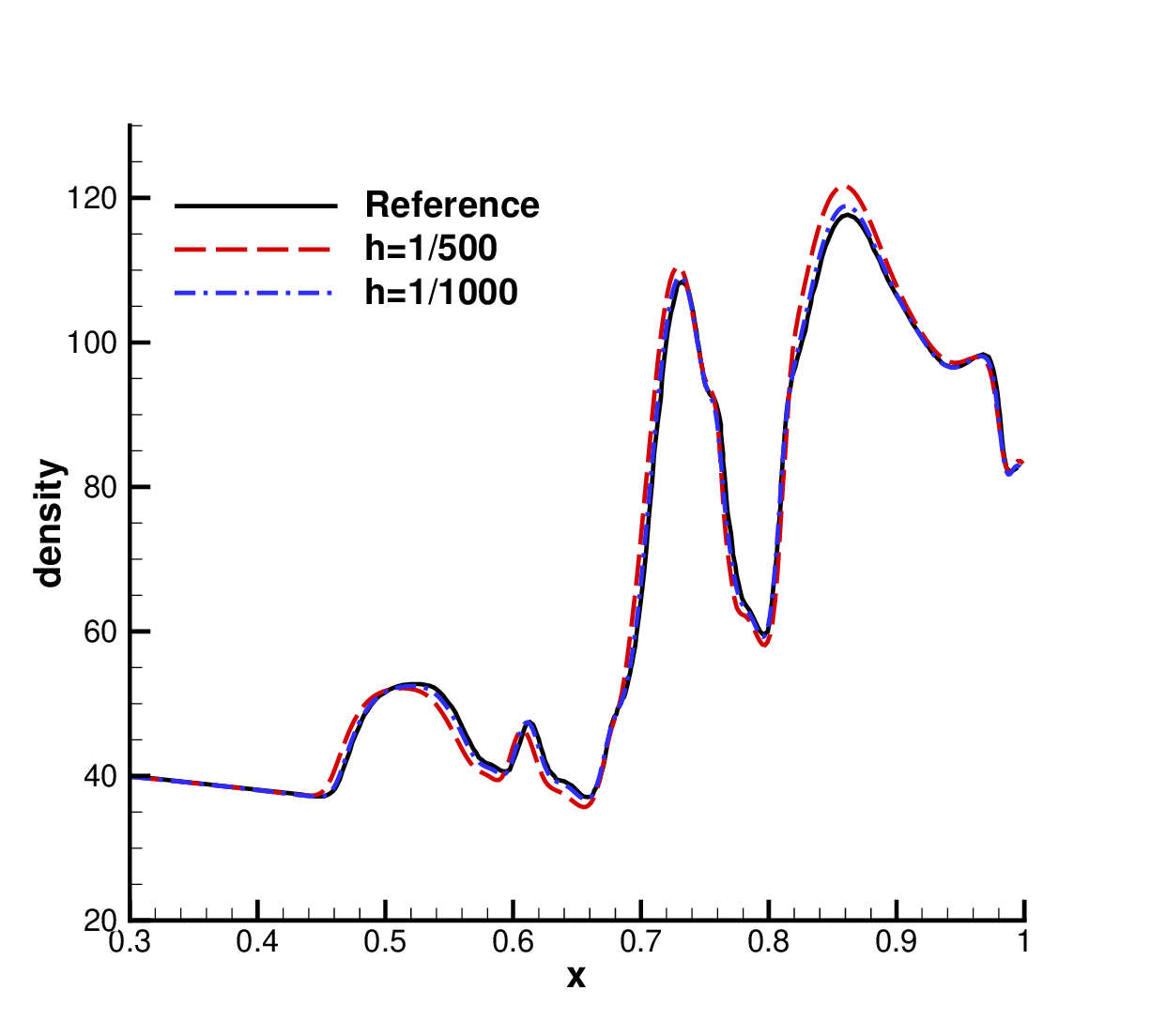}
	\caption{\footnotesize Reflected shock-boundary layer interaction: density along the bottom wall with the mesh size $h=1/500$ and $1/1000$ for $Re=200$. The reference data are from the HGKS method \cite{ZXL2018} with $h=1/1500$.}
	\label{fig:VSTP200den}
\end{figure}

The case of $\Re=200$ is tested first, where uniform meshes with $500\times250$ and $1000\times500$ cells are used. The density contours at $t=1.0$ are presented in \Cref{fig:VSTPRe200} with the two different mesh resolutions. It is observed that the complex flow structures are well resolved, including the lambda shock and the vortex structures. The density along the bottom wall is presented in \Cref{fig:VSTP200den}, which match very well with the reference data. \Cref{Tab:VSTP} presents the height of the primary vortex predicted by the current scheme, achieving a good agreement with the reference data obtained by a high-order GKS (HGKS) method in \cite{ZXL2018} with the mesh size $h=1/1500$.

\begin{figure}[!htb]
	\centering
	\subfigure{	\includegraphics[width=0.9\linewidth]{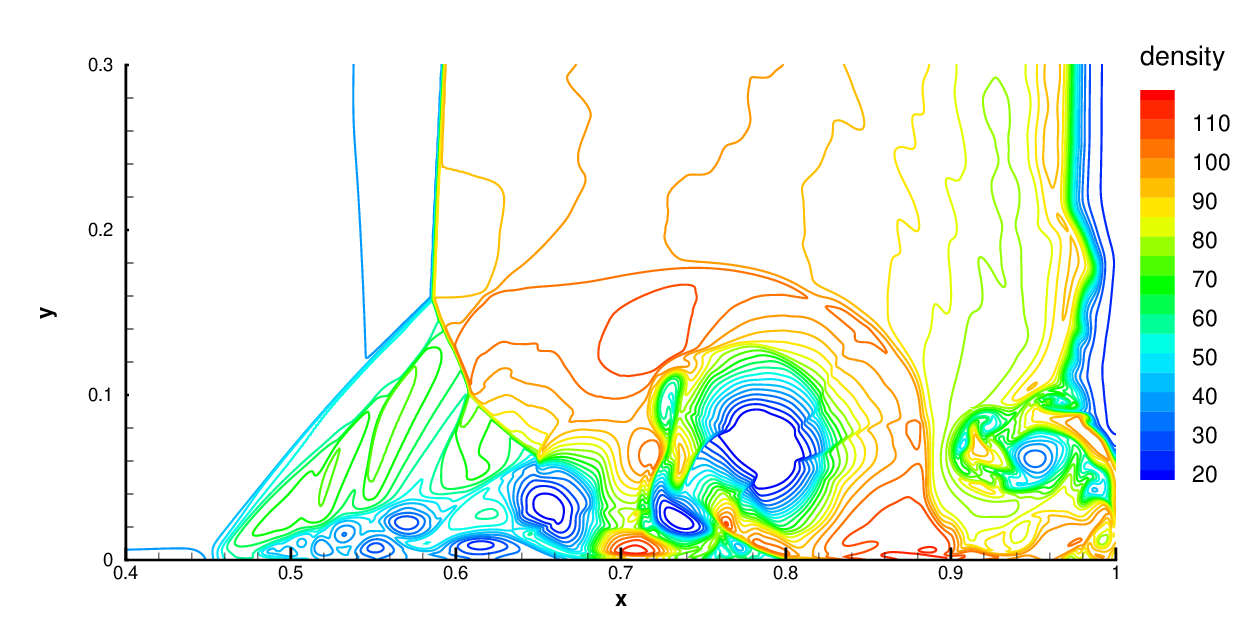}}
	\subfigure{	\includegraphics[width=0.9\linewidth]{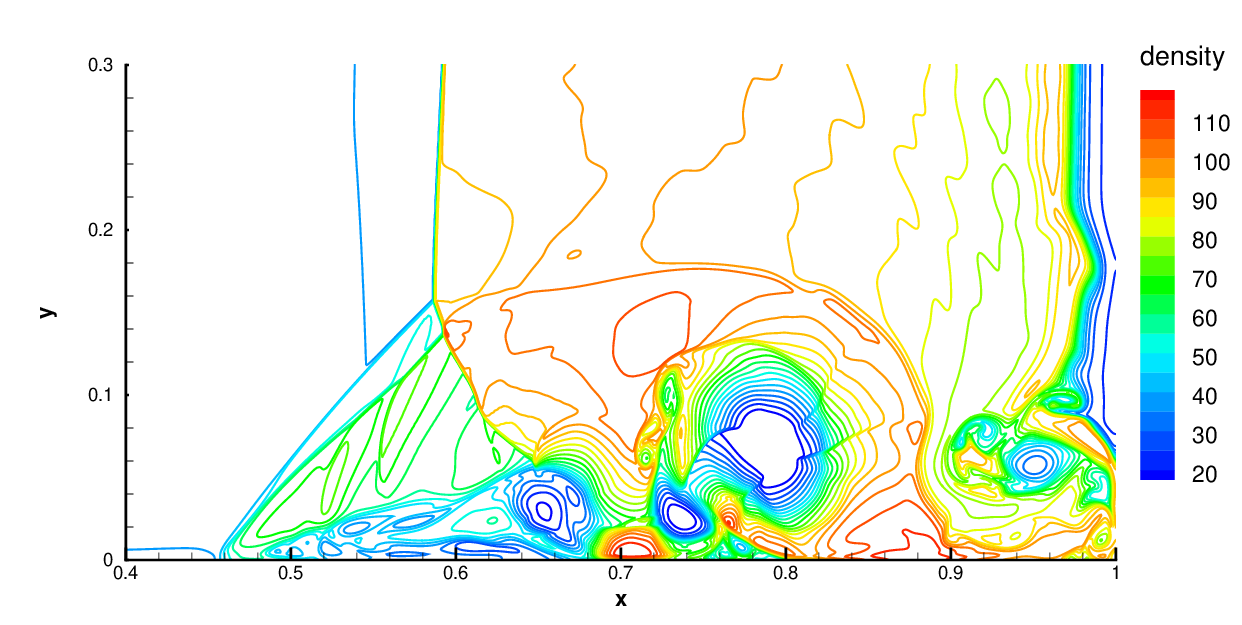}}
	\caption{\footnotesize Reflected shock-boundary layer interaction: density contours at $t=1.0$ for $\Re=1000$ with the mesh size $h=1/1500$ (top) and $1/2000$ (bottom). 20 uniform contours are displayed as the density ranges from 20 to 115.}
	\label{fig:VSTPRe1000}
\end{figure}

\begin{figure}[!htb]
	\centering
	\includegraphics[width=0.9\linewidth]{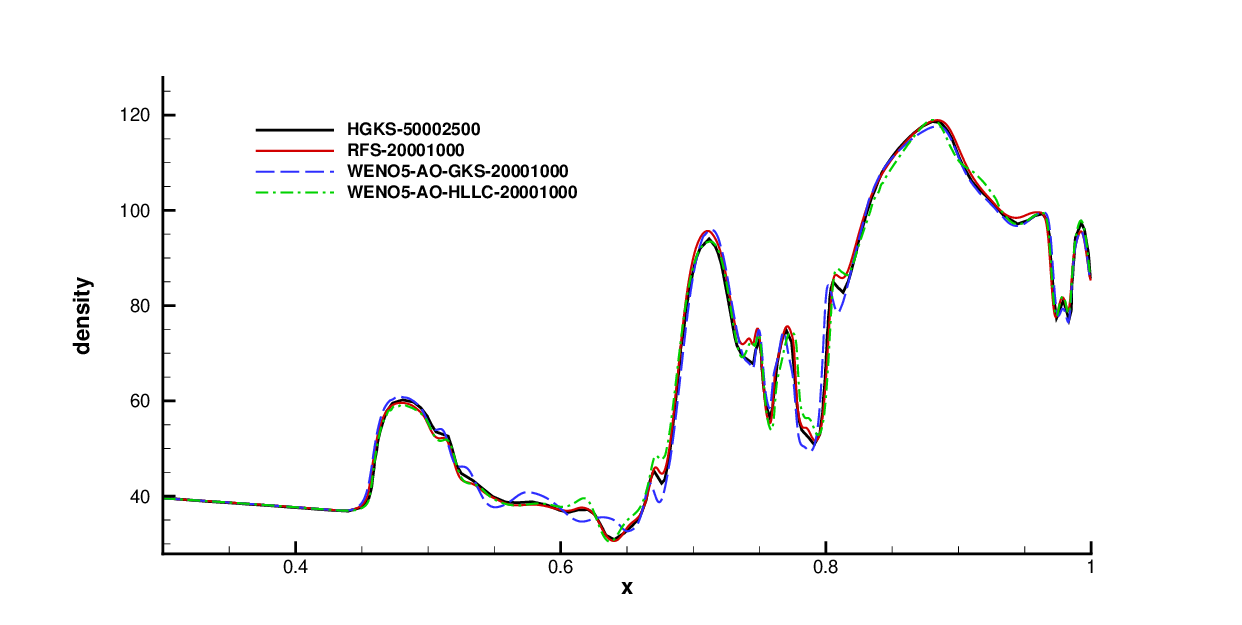}
	\includegraphics[width=0.9\linewidth]{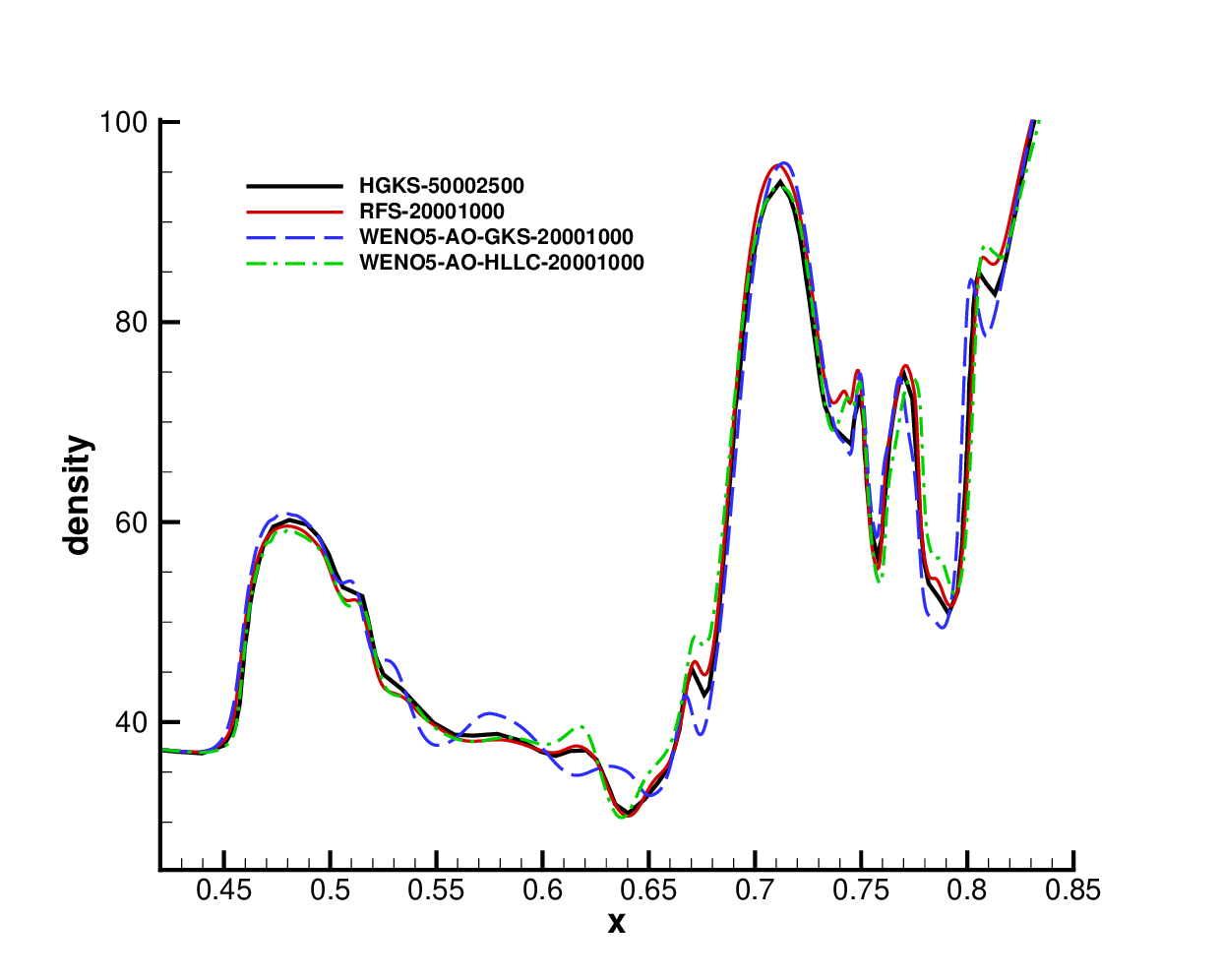}
	\caption{\footnotesize Reflected shock-boundary layer interaction: density along the bottom wall with $2000\times1000$ cells for $Re=1000$. The reference data are obtained from \cite{ZXL2018} and \cite{YJSX2022}.}
	\label{fig:VSTP1000den}
\end{figure}

For the case of $\Re=1000$, the flow structure becomes more complicated. Uniform meshes with $1500\times750$ and $2000\times1000$ cells are used. The density contours at $t=1.0$ from the current scheme are displayed in \Cref{fig:VSTPRe1000} with the mesh size $h=1/1500$ and $1/2000$, which match well with each other. The main flow structures can be captured very well by the current scheme. \Cref{fig:VSTP1000den} presents the density along the bottom wall and its local enlargement, obtained with $2000\times1000$ cells. The reference data are obtained by a HGKS method with $5000\times2500$ cells \cite{ZXL2018}, and fifth-order finite-volume methods based on the GKS solver and the HLLC solver with $2000\times1000$ cells \cite{YJSX2022}, respectively. The result obtained by the current scheme matches well with the reference data. 

Though the current scheme has only second-order accuracy, it resolves well such a complex flow, demonstrating the good performance of the current scheme in high-speed viscous flows.

\begin{table}[htbp]
	\centering
\caption{\footnotesize Reflected shock-boundary layer interaction: comparison of the height of the primary vortex for $\Re=200$. The reference data are from the HGKS method \cite{ZXL2018}.}
	\begin{tabular}{cccc}
		\toprule
		Scheme   
		&\multicolumn{2}{c}{RFS}   &  HGKS  \\
		\hline
		Mesh size  & $1/500$ & $1/1000$  &$1/1500$  \\
		\hline
		Height     & 0.163 & 0.165 &0.166 \\
		\bottomrule  
	\end{tabular}   
   \label{Tab:VSTP}  
\end{table}

\section{Discussions} \label{Discussions}

In this work, a GRP-based second-order finite volume scheme for the compressible Navier-Stokes equations is developed. 
The proposed scheme employs a hyperbolic relaxation model to compute the numerical flux. Both the convective and viscous fluxes in the governing PDEs are relaxed as stiff source terms in the relaxation model.

A key innovation of this method is that it avoids explicitly evolving relaxation variables in the finite volume discretization, using them only to construct numerical fluxes. This design significantly reduces computational cost and memory consumption. The GRP methodology is employed to provide both Riemann solutions and time derivatives. Thus, the effects of source terms are included in the development of the numerical flux. The resulting numerical procedure is termed the relaxation flux solver in the current paper.

This innovative use of the relaxation model allows the local determination of the parameters used in the model.
These parameters can vary dynamically at different cell interfaces and time steps, improving the adaptability to diverse flow regions and enhancing the robustness of the proposed scheme.

In the computation of time derivatives through the GRP methodology, the stiff source terms are treated implicitly, enabling the use of a standard hyperbolic CFL condition. This approach not only removes the severe parabolic time step restriction but also achieves second-order temporal accuracy through a single-stage evolution.

The robustness and accuracy of the proposed method were validated through extensive numerical tests covering a wide range of flow problems, from nearly incompressible regimes to supersonic flows with strong shocks, for both inviscid and viscous scenarios.

The present second-order accurate scheme has the potential to be extended to a high-order one in future studies, since recent advancements in LW type numerical methods have led to the development of high-order accurate finite volume schemes for compressible fluid flows \cite{LD2016,DL2018,PXLL2016,JZSX2020,LL2023}.
It is also worth pointing out that although the present scheme is designed for the compressible Navier-Stokes equations, it may be extended to other convection-diffusion equations, e.g., the resistive magneto-hydrodynamic model.

\appendix

\section{The homogeneity tensor} \label{tensor}
The matrices $\bfB_{11}(\bfu)$ and  $\bfB_{22}(\bfu)$ in \eqref{matricesB} are given by
\begin{equation}\notag
	\bfB_{11}(\bfu)=\frac{\mu}{\rho}
	\begin{pmatrix}
		0\quad\quad&0\quad\quad&0\quad\quad&0\quad \\
		-\frac{4}{3} u\quad\quad&\frac{4}{3}\quad\quad&0\quad\quad&0\quad \\
		- v\quad\quad&0\quad\quad&1\quad\quad&0\quad \\
		B_{4111}\quad\quad&(\frac{4}{3}-\frac{\gamma}{\Pr})u\quad\quad
		&(1-\frac{\gamma}{\Pr})v\quad\quad&\frac{\gamma}{\Pr}\quad
	\end{pmatrix}
\end{equation}
and
\begin{equation}\notag
	\bfB_{22}(\bfu)=\frac{\mu}{\rho}
	\begin{pmatrix}
		0\quad\quad&0\quad\quad&0\quad\quad&0\quad \\
		- u\quad\quad&1\quad\quad&0\quad\quad&0\quad \\
		-\frac{4}{3} v\quad\quad&0\quad\quad&\frac{4}{3}\quad\quad&0\quad \\
		B_{4212}\quad\quad&(1-\frac{\gamma}{\Pr})u\quad\quad
		&(\frac{4}{3}-\frac{\gamma}{\Pr})v\quad\quad&\frac{\gamma}{\Pr}\quad
	\end{pmatrix}
\end{equation}
where
\begin{equation*}
\begin{aligned}
	&B_{4111}=-\frac{4}{3} u^2- v^2-\frac{\gamma}{\Pr}(E-(u^2+v^2)),\\
	&B_{4212}=- u^2-\frac{4}{3} v^2-\frac{\gamma}{\Pr}(E-(u^2+v^2)).
\end{aligned}
\end{equation*}
The matrices $\bfB_{12}(\bfu)$ and  $\bfB_{21}(\bfu)$ are
\begin{equation*}
	\bfB_{12}(\bfu)=\frac{\mu}{\rho}
	\begin{pmatrix}
		0\quad\quad&0\quad\quad&0\quad\quad&0\quad \\
		\frac{2}{3} v\quad\quad&0\quad\quad&-\frac{2}{3}\quad\quad&0\quad \\
		- u\quad\quad&1\quad\quad&0\quad\quad&0\quad \\
		-\frac{1}{3} uv\quad\quad& v\quad\quad&-\frac{2}{3} u\quad\quad&0\quad 
	\end{pmatrix},
\end{equation*}
and
\begin{equation*}
	\bfB_{21}(\bfu)=\frac{\mu}{\rho}
	\begin{pmatrix}
		0\quad\quad&0\quad\quad&0\quad\quad&0\quad \\
		- v\quad\quad&0\quad\quad&1\quad\quad&0\quad \\
		\frac{2}{3} u\quad\quad&-\frac{2}{3}\quad\quad&0\quad\quad&0\quad \\
		-\frac{1}{3} uv\quad\quad&-\frac{2}{3} v\quad\quad& u\quad\quad&0\quad 
	\end{pmatrix}.
\end{equation*}

\section{Space data reconstruction} \label{datarecon} 
This appendix is devoted to the space data reconstruction procedure for the finite volume scheme \eqref{FV}. A second-order accurate space data reconstruction similar to the one proposed in \cite{BLW2006} is adopted. Here, we apply the linear reconstruction to primitive variables $\bfQ=(\rho,u,v,T)^\top$.

At time level $t=t^{n+1}$, $\bfQ$ is linearly reconstructed in each cell $\Omega_{i,j}$ as
\begin{equation}
	\bfQ_{i,j}^{n+1}(x,y)
	=\bftQ_{i,j}^{n+1}+\bigg(\frac{\pt\bfQ}{\pt x }\bigg)_{i,j}^{n+1}(x-x_i)
	+\bigg(\frac{\pt\bfQ}{\pt y }\bigg)_{i,j}^{n+1}(y-y_j),
	\quad (x,y)\in \Omega_{i,j}.  \label{reQ}
\end{equation}
Here, the cell centered value $\bftQ_{i,j}^{n+1}$ is computed by using the cell average values $\bfbu_{i,j}^{n+1}$ obtained in \eqref{FV}, i.e.,
\begin{equation}\notag
	\bftQ_{i,j}^{n+1}=\bfQ(\bfbu_{i,j}^{n+1}).
\end{equation}

For smooth flow problems, the slope $\big(\frac{\pt\bfQ}{\pt x }\big)_{i,j}^{n+1}$ is linearly approximated by
\begin{equation}
	\bigg(\frac{\pt\bfQ}{\pt x }\bigg)_{i,j}^{n+1}
    =\frac{\bftQ_{i+1,j}^{n+1}-\bftQ_{i-1,j}^{n+1}}{2\Delta x_{i}}.	\label{recon}
\end{equation} 
In the presence of discontinuities, we employ the minmod limiter \cite{BLW2006,CS1989,vanLeer1979} to construct the slope $\big(\frac{\pt\bfQ}{\pt x }\big)_{i,j}^{n+1}$ as
\begin{equation}
	\bigg(\frac{\pt\bfQ}{\pt x }\bigg)_{i,j}^{n+1}
=\mathrm{minmod}\bigg(\frac{\alpha(\bftQ^{n+1}_{i+1,j}-\bftQ^{n+1}_{i,j})}{\frac{1}{2}(\Delta x_{i+1}+\Delta x_{i})},
	\frac{\bfQ^{n+1,-}_{i+\hf,j}-\bfQ^{n+1,-}_{i-\hf,j}}
    {\Delta x_i},
    \frac{\alpha(\bftQ^{n+1}_{i,j}-\bftQ^{n+1}_{i-1,j})}{\frac{1}{2}(\Delta x_i+\Delta x_{i-1})}\bigg),
	\label{minmod}
\end{equation} 
where $\bfQ^{n+1,-}_{i+\hf,j}$ and $\bfQ^{n+1,-}_{i+\hf,j}$ are already obtained through the GRP methodology in \eqref{gradient1-1}, and the coefficient $\alpha\in[0,2)$ is a user-tuned parameter.

The slope in the $y$-direction, $\big(\frac{\pt\bfQ}{\pt y }\big)_{i,j}^{n+1}$, can be computed similarly. By transforming the reconstructed piecewise linear initial data of $\bfQ$ in \eqref{reQ} into those of $\bfu$, the piecewise smooth initial data $\bfu_{i+1,j}^{n+1}(x,y)$ and $\bfu_{i,j}^{n+1}(x,y)$ to be used in the IVP \eqref{FVGRP} are obtained.

\section{Explicit components of numerical fluxes} \label{terms}
The explicit parts of numerical fluxes in \eqref{eq:impl-evol-u} are given below. For terms at $t^n$, we have
\begin{equation}\label{eq:def-R1u}
    \begin{aligned}
        R^{u}_{i,j,1}=&-\frac{(1-\omega^n_{i+\hf,j})(\bfv_2)_{i+\hf,j}^{n,*}
            -(1-\omega^n_{i-\hf,j})(\bfv_2)_{i-\hf,j}^{n,*}}{\Delta x_i}\\
            &-\frac{(1-\omega^n_{i,j+\hf})(\bfw_2)_{i+\hf,j}^{n,*}-
            (1-\omega^n_{i,j-\hf})(\bfw_2)_{i-\hf,j}^{n,*}}{\Delta y_j},
    \end{aligned}
\end{equation}
and
\begin{equation}\label{eq:def-R23u}
    \begin{aligned}
         R^{u}_{i,j,2}=&\frac{(1-\omega^n_{i+\hf,j})(a^{n}_{i+\hf,j})^2\Delta t
        (\frac{\pt\bfu_2}{\pt x})_{i+\hf,j}^{n,*}
            -(1-\omega^n_{i-\hf,j})(a^{n}_{i-\hf,j})^2\Delta t
        (\frac{\pt\bfu_2}{\pt x})_{i-\hf,j}^{n,*}}{\Delta x_i}\\
          &+\frac{(1-\omega^n_{i,j+\hf})(a^{n}_{i,j+\hf})^2\Delta t
        (\frac{\pt\bfu_2}{\pt x})_{i,j+\hf}^{n,*}
            -(1-\omega^n_{i,j-\hf})(a^{n}_{i,j-\hf})^2\Delta t
        (\frac{\pt\bfu_2}{\pt x})_{i,j-\hf}^{n,*}}{\Delta y_j},\\
         R^{u}_{i,j,3}=&-\frac{1}{2}\frac{\omega^n_{i+\hf,j}(\bfH_2)_{i+\hf,j}^{n}-\omega^n_{i-\hf,j}(\bfH_2)_{i-\hf,j}^{n}}{\Delta x_i}\\
&-\frac{1}{2}\frac{\omega^n_{i,j+\hf}(\mathbf{K}_2)_{i,j+\hf}^{n}-\omega^n_{i,j-\hf}(\mathbf{K}_2)_{i,j-\hf}^{n}}{\Delta y_j},
    \end{aligned}
\end{equation}
where $\mathbf{K}=\bfg_c-\bfg_v$, $(\cdot)_2$ is the second component of the vector, and
\begin{equation*}
		\omega^n_{i+\hf,j} = \frac{\Delta t}{2\ep^n_{i+\hf,j}+\Delta t },
		\quad \omega^n_{i,j+\hf}  = \frac{\Delta t}{2\ep^n_{i,j+\hf}+\Delta t}.
\end{equation*}
For terms at $t^{n+1}$, we have
\begin{equation}\label{eq:def-R456u}
    \begin{aligned}
   R^{u}_{i,j,4}=&-\frac{1}{2}\frac{\omega^n_{i+\hf,j}((\bff_c)_2)_{i+\hf,j}^{n+1,-}-\omega^n_{i-\hf,j}((\bff_c)_2)_{i-\hf,j}^{n+1,-}}{\Delta x_i}\\
&-\frac{1}{2}\frac{\omega^n_{i,j+\hf}((\bfg_c)_2)_{i,j+\hf}^{n+1,-}-\omega^n_{i,j-\hf}((\bfg_c)_2)_{i,j-\hf}^{n+1,-}}{\Delta y_j},\\
 R^{u}_{i,j,5}=&-\frac{1}{3}
\frac{\mu(T^{n+1,-}_{i+\hf,j})\omega^n_{i,j+\hf}(\frac{\pt v}{\pt y})^{n+1,-}_{i+\hf,j}-
\mu(T^{n+1,-}_{i-\hf,j})\omega^n_{i-\hf,j}(\frac{\pt v}{\pt y})^{n+1,-}_{i-\hf,j}}{\Delta x_i}\\
&+\frac{1}{2}\frac{\mu(T^{n+1,-}_{i,j+\hf})\omega^n_{i,j+\hf}(\frac{\pt v}{\pt x})_{i,j+\hf}^{n+1,-}
-\mu(T^{n+1,-}_{i,j-\hf})\omega^n_{i,j-\hf}(\frac{\pt v}{\pt x})_{i,j-\hf}^{n+1,-}}{\Delta y_j},\\
R^{u}_{i,j,6}=&\frac{1}{3}
\frac{\mu(T^{n+1,-}_{i+\hf,j})\omega^n_{i+\hf,j}R^{u}_{i,j,61}
-\mu(T^{n+1,-}_{i-\hf,j})\omega^n_{i-\hf,j}R^{u}_{i,j,62}}{\Delta x_i}\\
&+\frac{1}{4}
\frac{\mu(T^{n+1,-}_{i,j+\hf})\omega^n_{i,j+\hf}R^{u}_{i,j,63}
-\mu(T^{n+1,-}_{i,j-\hf})\omega^n_{i,j-\hf}R^{u}_{i,j,64}}{\Delta y_j},\\
    \end{aligned}
\end{equation}
with
\begin{equation}\notag
    \begin{aligned}
R^{u}_{i,j,61}=&\big(1-\frac{\Delta x_i}{2\Delta x_{i+\hf}}\big)
\bigg(\frac{\pt u}{\pt x}\bigg)_{i+\hf-0,j}^{n+1,-}
+\big(1-\frac{\Delta x_{i+1}}{2\Delta x_{i+\hf}}\big)
\bigg(\frac{\pt u}{\pt x}\bigg)_{i+\hf+0,j}^{n+1,-},\\
R^{u}_{i,j,62}=&\big(1-\frac{\Delta x_{i-1}}{2\Delta x_{i-\hf}}\big)
\bigg(\frac{\pt u}{\pt x}\bigg)_{i-\hf-0,j}^{n+1,-}
+\big(1-\frac{\Delta x_{i}}{2\Delta x_{i-\hf}}\big)
\bigg(\frac{\pt u}{\pt x}\bigg)_{i-\hf+0,j}^{n+1,-},\\
R^{u}_{i,j,63}=&\big(1-\frac{\Delta y_j}{2\Delta y_{j+\hf}}\big)
\bigg(\frac{\pt u}{\pt y}\bigg)_{i,j+\hf-0}^{n+1,-}
+\big(1-\frac{\Delta y_{j+1}}{2\Delta y_{j+\hf}}\big)
\bigg(\frac{\pt u}{\pt y}\bigg)_{i,j+\hf+0}^{n+1,-},\\
R^{u}_{i,j,64}=&\big(1-\frac{\Delta y_{j-1}}{2\Delta y_{j-\hf}}\big)
\bigg(\frac{\pt u}{\pt y}\bigg)_{i,j-\hf-0}^{n+1,-}
+\big(1-\frac{\Delta y_{j}}{2\Delta y_{j-\hf}}\big)
\bigg(\frac{\pt u}{\pt y}\bigg)_{i,j-\hf+0}^{n+1,-}.
    \end{aligned}
\end{equation}

By symmetry, the coefficient matrix $\fM_{v}$ and the right-hand vector $\vb_v$ have a similar form of $\fM_{u}$ and $\vb_u$, respectively, and the details are omitted here.

After solving the linear equation system of $u$ and $v$, the cell interface values $\big(\frac{\pt u}{\pt x}\big)^{n+1,-}_{i+\hf,j}$ and $\big(\frac{\pt v}{\pt y}\big)^{n+1,-}_{i,j+\hf}$ are obtained by \eqref{gradient1-1} and \eqref{implicit0}. Thus, for the numerical flux of the energy conservation equation, the viscosity effect terms are already obtained and can be regarded as explicit components. By the same manipulation, the coefficient matrix of $T$ is
\begin{equation} \label{Jacobian-T}
	(\fM_T)_{k,l}=\begin{cases}
		\kappa(T^{n+1,-}_{i-\hf,j})\frac{\Delta t}{2\Delta x_i \Delta x_{i-\hf}}\omega^n_{i-\hf,j},
		\quad\quad &\text{if}\,\, l=k-I,\\[2mm]
		\kappa(T^{n+1,-}_{i+\hf,j})\frac{\Delta t}{2\Delta x_i \Delta x_{i+\hf}}\omega^n_{i+\hf,j},
		\quad\quad &\text{if}\,\, l=k+I,\\[2mm]  
		\kappa(T^{n+1,-}_{i,j-\hf})\frac{\Delta t}{2\Delta y_j \Delta y_{j-\hf}}\omega^n_{i,j-\hf},
		\quad\quad &\text{if}\,\, l=k-1,\\[2mm]
		\kappa(T^{n+1,-}_{i,j+\hf})\frac{\Delta t}{2\Delta y_j \Delta y_{j+\hf}}\omega^n_{i,j+\hf},
		\quad\quad &\text{if}\,\,l=k+1,\\[2mm] 
		\frac{\bar{\rho}^{n+1}_{i,j}}{\gamma-1}
		+(\fM_T)^k_{k-I}+(\fM_T)^k_{k+I}\\
		+(\fM_T)^k_{k-1}+(\fM_T)^k_{k+1},     &\text{if}\,\, l=k,\\[2mm]
		0, \quad\quad\quad\quad          &\text{otherwise}.
	\end{cases}	         
\end{equation} 

The $k$-th entry of the right-hand side vector $\vb_T$ is
\begin{equation}
(\vb_T)_k=\frac{\bar\rho_{i,j}^{n}\tilde{T}_{i,j}^{n}}{\gamma-1}+\Delta t\sum_{s=1}^{6}R^{T}_{i,j,s},
\end{equation}
where the cell index $(i,j)$ satisfies $(i-1)\times I+j=k$. At the time level $t^n$,
\begin{equation}\label{eq:def-R123T}
    \begin{aligned}
        R^{T}_{i,j,1}=&-\frac{(1-\omega^n_{i+\hf,j})(\bfv_4)_{i+\hf,j}^{n,*}
            -(1-\omega^n_{i-\hf,j})(\bfv_4)_{i-\hf,j}^{n,*}}{\Delta x_i}\\
            &-\frac{(1-\omega^n_{i,j+\hf})(\bfw_4)_{i+\hf,j}^{n,*}-
            (1-\omega^n_{i,j-\hf})(\bfw_4)_{i-\hf,j}^{n,*}}{\Delta y_j},\\
            R^{T}_{i,j,2}=&\frac{(1-\omega^n_{i+\hf,j})(a^{n}_{i+\hf,j})^2\Delta t
        (\frac{\pt\bfu_4}{\pt x})_{i+\hf,j}^{n,*}
            -(1-\omega^n_{i-\hf,j})(a^{n}_{i-\hf,j})^2\Delta t
        (\frac{\pt\bfu_4}{\pt x})_{i-\hf,j}^{n,*}}{\Delta x_i}\\
          &+\frac{(1-\omega^n_{i,j+\hf})(a^{n}_{i,j+\hf})^2\Delta t
        (\frac{\pt\bfu_4}{\pt x})_{i,j+\hf}^{n,*}
            -(1-\omega^n_{i,j-\hf})(a^{n}_{i,j-\hf})^2\Delta t
        (\frac{\pt\bfu_4}{\pt x})_{i,j-\hf}^{n,*}}{\Delta y_j},\\
         R^{T}_{i,j,3}=&-\frac{1}{2}\frac{\omega^n_{i+\hf,j}(\bfH_4)_{i+\hf,j}^{n}-\omega^n_{i-\hf,j}(\bfH_4)_{i-\hf,j}^{n}}{\Delta x_i}
-\frac{1}{2}\frac{\omega^n_{i,j+\hf}(\mathbf{K}_4)_{i,j+\hf}^{n}-\omega^n_{i,j-\hf}(\mathbf{K}_4)_{i,j-\hf}^{n}}{\Delta y_j}.
    \end{aligned}
\end{equation}
At the time level $t^{n+1}$,
\begin{equation}\label{eq:def-R456T}
    \begin{aligned}
   R^{T}_{i,j,4}=&-\frac{1}{2}\frac{\omega^n_{i+\hf,j}((\bff_c)_4)_{i+\hf,j}^{n+1,-}-\omega^n_{i-\hf,j}((\bff_c)_4)_{i-\hf,j}^{n+1,-}}{\Delta x_i}\\
&-\frac{1}{2}\frac{\omega^n_{i,j+\hf}((\bfg_c)_4)_{i,j+\hf}^{n+1,-}-\omega^n_{i,j-\hf}((\bfg_c)_4)_{i,j-\hf}^{n+1,-}}{\Delta y_j},\\
R^{T}_{i,j,5}=&\frac{1}{2}
\frac{\omega^n_{i,j+\hf}\big(u(\bff_v)_2+v(\bff_v)_3\big)^{n+1,-}_{i+\hf,j}
-\omega^n_{i-\hf,j}\big(u(\bff_v)_2+v(\bff)_3\big)^{n+1,-}_{i-\hf,j}}{\Delta x_i}\\
&+\frac{1}{2}
\frac{\omega^n_{i,j+\hf}\big(u(\bfg_v)_2+v(\bfg_v)_3\big)_{i,j+\hf}^{n+1,-}
-\omega^n_{i,j-\hf}\big(u(\bfg_v)_2+v(\bfg_v)_3\big)_{i,j-\hf}^{n+1,-}}{\Delta y_j},\\
R^{T}_{i,j,6}=&\frac{1}{4}
\frac{\kappa(T^{n+1,-}_{i+\hf,j})\omega^n_{i+\hf,j}R^{T}_{i,j,61}
-\kappa(T^{n+1,-}_{i-\hf,j})\omega^n_{i-\hf,j}R^{T}_{i,j,62}}{\Delta x_i}\\
&+\frac{1}{4}
\frac{\kappa(T^{n+1,-}_{i,j+\hf})\omega^n_{i,j+\hf}R^{T}_{i,j,63}
-\kappa(T^{n+1,-}_{i,j-\hf})\omega^n_{i,j-\hf}R^{T}_{i,j,64}}{\Delta y_j},\\
    \end{aligned}
\end{equation}
with
\begin{equation}\notag
    \begin{aligned}
R^{T}_{i,j,61}=&\big(1-\frac{\Delta x_i}{2\Delta x_{i+\hf}}\big)
\bigg(\frac{\pt T}{\pt x}\bigg)_{i+\hf-0,j}^{n+1,-}
+\big(1-\frac{\Delta x_{i+1}}{2\Delta x_{i+\hf}}\big)
\bigg(\frac{\pt T}{\pt x}\bigg)_{i+\hf+0,j}^{n+1,-},\\
R^{T}_{i,j,62}=&\big(1-\frac{\Delta x_{i-1}}{2\Delta x_{i-\hf}}\big)
\bigg(\frac{\pt T}{\pt x}\bigg)_{i-\hf-0,j}^{n+1,-}
+\big(1-\frac{\Delta x_{i}}{2\Delta x_{i-\hf}}\big)
\bigg(\frac{\pt T}{\pt x}\bigg)_{i-\hf+0,j}^{n+1,-},\\
R^{T}_{i,j,63}=&\big(1-\frac{\Delta y_j}{2\Delta y_{j+\hf}}\big)
\bigg(\frac{\pt T}{\pt y}\bigg)_{i,j+\hf-0}^{n+1,-}
+\big(1-\frac{\Delta y_{j+1}}{2\Delta y_{j+\hf}}\big)
\bigg(\frac{\pt T}{\pt y}\bigg)_{i,j+\hf+0}^{n+1,-},\\
R^{T}_{i,j,64}=&\big(1-\frac{\Delta y_{j-1}}{2\Delta y_{j-\hf}}\big)
\bigg(\frac{\pt T}{\pt y}\bigg)_{i,j-\hf-0}^{n+1,-}
+\big(1-\frac{\Delta y_{j}}{2\Delta y_{j-\hf}}\big)
\bigg(\frac{\pt T}{\pt y}\bigg)_{i,j-\hf+0}^{n+1,-}.\\
    \end{aligned}
\end{equation}

\vspace{2mm}
\section*{Acknowledgement}
This research is supported by NSFC (Nos. 12371391, 12031001).

\vspace{2mm}

\bibliographystyle{siamplain}

\end{document}